\documentclass[11pt]{article}

\input{amssym.def}
\input{amssym}
\usepackage{eufrak}

\newtheorem{theorem}{Theorem}

\newtheorem{conjecture}[theorem]{Conjecture}
\newtheorem{corollary}[theorem]{Corollary}

\newtheorem{definition}[theorem]{Definition}
\newtheorem{example}[theorem]{Example}

\newtheorem{lemma}[theorem]{Lemma}

\newtheorem{proposition}[theorem]{Proposition}
\newtheorem{remark}[theorem]{Remark}

\newlength{\coefkappa}
\begin{document}

\def\h{{\hbar}}
\def\Cas{{\bf Cas}}
\def\Op{{\rm Op}}
\def\Caslm{{\bf Cas}_{(\la,m)}}
\def\ch{{\cal CH}}
\def\chc{{\cal CH^{\chi}}}
\def\lhq{\ifmmode {\cal L}_{\hbar,q}\else ${\cal L}_{\hbar,q}$\fi}
\def\lq{{\cal L}_{q}}
\def\loq{{\cal L}_{1,q}}
\def\lqh{\ifmmode {\cal L}_{\hbar,q}\else ${\cal L}_{\hbar,q}$\fi}
\def\slhq{{\cal SL}_{\h, q}}
\def\lc{{\cal L}_{\h, q}^{\chi}}
\def\gg{\mbox{$\frak g$}}

\def\gh{\gg_{\h}}
\def\Ih{I_{\h}}
\def\lhc{{\cal L}_{\h}^{\chi}}
\def\diag{\rm diag}
\def\li{{\cal L}_i}
\def\lm{{ L}_{(m)}}
\def\llm{{ L}_{(\la,m)}}
\def\ppm{p_{(m)}}
\def\Id{I}
\def\Inv{\rm Inv}
\def\vl{V_{\la}}
\def\rl{\rho_{\la}}
\def\rlt{\rho_{\la}^{(2)}}
\def\mult{\rm mult}
\def\MM{\cal M}
\def\Min{rm Min}
\def\nl{n_\la}
\def\Fl{\rm Fl}
\def\Rep{\rm Rep}
\def\zl{Z(\lhq)}
\def\lhqc{{\cal L}_{\h,q}^{\chi}}
\def\KK{Q}
\def\vm{V_{(m)}}
\def\Lm{L_{(m)}}
\def\lmod{{\lhq}\mbox{--}\Mod}
\def\A{\cal A}
\def\k{\mathbf{k}}
\def\qq{q^{-1}}
\def\lq{{\cal L}_q}
\def\lqc{{\cal L}_q^{\chi}}
\def\chm{{\cal CH}_{(m)}}
\def\ola{\overline{\lambda}}
\def\Zq{\Bbb{Z}_q}
\def\vmd{{\frak D}}

\def\sll{s_{\la}}
\def\glh{gl(2)_{\h}}
\def\gq{\gggg_{q}}
\def\gqq{\gggg_{q}^{\ot 2}}
\def\Ugl{U(gl(V))}
\def\Usl{U(sl(V))}
\def\Ush{U(sl(2)_{\h})}
\def\Zgh{Z[U(gl(2)_{\h})]}
\def\Zsh{Z(U(sl(2)_{\h}))}
\def\uq{U_q(sl(2))}
\def\Uq{U_q(sl(n))}
\def\ugh{U(\gggg_{\h})}
\def\uqs{U_q(sl(n))}
\def\ug{U(\gggg)}
\def\C{{\Bbb C}}
\def\B{{\cal C}}
\def\eh{e(\h)}
\def\mh{M(\h)}
\def\bl{{B_{\la}}}
\def\cl{{C_{\la}}}
\def\cp{{\Bbb CP}}
\def\N{{\Bbb N}}
\def\R{{\Bbb R}}
\def\K{{\Bbb K}}
\def\Z{{\Bbb Z}}
\def\O{{\cal O}}
\def\A{{\cal A}}
\def\Ah{{\cal A}_\hbar}
\def\kO{{\K(\cal O)}}
\def\De{\Delta}
\def\de{\delta}
\def\la{\lambda}
\def\La{\Lambda}
\def\co{{\cal O}}

\def\Pp{P_+}
\def\Pm{P_-}
\def\Ppm{P_{\pm}}

\def\lp{\Lambda_+(V)}
\def\lmm{\Lambda_-(V)}
\def\lpm{\Lambda_{\pm}(V)}
\def\lpmk{\Lambda_{\pm}^k(V)}
\def\lmk{\Lambda_{-}^k(V)}
\def\lmpp{\Lambda_{-}^{p-1}(V)}
\def\lmp{\Lambda_{-}^p(V)}

\def\l{{\bf l}}
\def\opi{\overline \pi}
\def\lpml{\Lambda_{\pm}^l(V)}
\def\lmp{\Lambda_{-}^p(V)}
\def\qm{q^{-1}}
\def\ve{\varepsilon}
\def\ep{\epsilon}
\def\al{{\alpha}}
\def\ah{{ A}_{\h}}
\def\aq{{ A}_{q}}
\def\ahq{{ A}_{\h,q}}
\def\ot{\otimes}
\def\om{\omega}
\def\Om{{\cal O}_M}
\def\ii{\rm i}
\def\sr{S^2_r}
\def\Ml{M^{\lambda}}
\def\End{{\rm End\, }}
\def\Mat{{\rm Mat }}
\def\Mor{{\rm Mor\, }}
\def\Trr{{\bf Tr}}
\def\Tr{{\rm Tr}}
\def\Ind{{\rm Ind\, }}
\def\prePic{{\rm prePic }}
\def\Mod{{\rm Mod}}
\def\Uqmod{U_q(sl(n))-{\rm Mod}}
\def\uqmod{U_q(sl(2))-{\rm Mod}}
\def\slm{sl(n)-{\rm Mod}}
\def\Hom{{\rm End\, }}
\def\Sym{{\rm Sym\, }}
\def\Ker{{\rm Ker\, }}
\def\Im{{\rm Im\, }}
\def\Ob{{\rm Ob\, }}

\def\id{{\rm id\, }}
\def\dim{{\rm dim}}
\def\ddim{{\rm dim}}
\def\span{{\rm span\, }}
\def\ttr{{\bf tr\, }}
\def\ad{{\rm ad\, }}
\def\cotr{{\rm cotr\, }}
\def\rk{{\rm rk\, }}
\def\Vect{{\rm Vect\,}}
\def\Span{{\rm Span}}
\def\Fun{{\rm Fun\,}}
\def\vv{V^{\ot 2}}
\def\vvv{V^{\ot 3}}
\def\vl{V_{\la}}
\def\lla{\bf{\la}}

\def\bea{\begin{eqnarray}}
\def\eea{\end{eqnarray}}
\def\be{\begin{equation}}
\def\ee{\end{equation}}
\def\nn{\nonumber}

\makeatletter
\renewcommand{\theequation}{{\thesection}.{\arabic{equation}}}
\@addtoreset{equation}{section} \makeatother

\title{\hfill\parbox{0.23\hsize}{
\normalsize\it In memory of our \\
friend J.Donin}\\
\bigskip
Geometry of non-commutative orbits\\
related to Hecke symmetries}
\author{D.~Gurevich\\
{\small\it ISTV, Universit\'e de Valenciennes, 59304 Valenciennes,
France}\\
\rule{0pt}{7mm}P.~Saponov,\\
{\small\it Theory Department of Institute for High Energy Physics,
142281 Protvino, Russia}}

\maketitle

\begin{abstract}
To some braiding $R$ of Hecke type (a Hecke symmetry) we put into
correspondence an associative algebra called the modified Reflection
Equation Algebra (mREA). We construct a series of matrices $L_{(m)},
\,m=1,2,...$ with entries belonging to such an algebra so that
each of them satisfies a version of the Cayley-Hamilton identity
with central coefficients.

We also consider some quotients of the mREA which are called
the non-commutative orbits. For each of these orbits we construct
a large family of projective modules. In such a family we introduce
an algebraic structure which is close to that of $K^0(\Fl(\C^n))$.
This algebraic structure respects an equivalence relation motivated
by a "quantum" trace compatible with the initial Hecke symmetry $R$.
For a subclass of non-commutative orbits we compute the spectrum of
central elements $\Tr_R\,L_{(m)}^k,\, k\in {\Bbb N}$ of the mREA
\end{abstract}

{\bf AMS Mathematics Subject Classification, 1991:} 17B37, 81R50


\section{Introduction}
\label{sec:int}

As was shown in \cite{Se}, the category of locally trivial vector
bundles over an affine regular algebraic variety $X$ is equivalent
to that of finitely generated projective $\K(X)$--modules where
$\K(X)$ is the coordinate algebra of $X$. (In what follows the
ground field $\K$ is assumed to be $\C$ or $\R$, the latter case
will be specified each time.) A similar equivalence is valid for
compact smooth varieties (cf. \cite{Ro}). This is the reason why
in the non-commutative geometry finitely generated projective 
modules over non-commutative algebras are considered as appropriate
analogs of vector bundles. Thus, the problem of constructing and
classifying such $\A$-modules for a given non-commutative (NC)
algebra $\A$ is of great interest. Unfortunately, besides
$\C^*$-algebras very few examples of algebras with a significant
family of projective modules are known.

In \cite{GS1} we have suggested a way of constructing projective
$\A$-modules via an NC version of the Cayley-Hamilton (CH)
identities. The existence of these identities  (as well as the
Newton identities) is the very remarkable property of the
so-called reflection equation algebras (REA) or their modified
versions (mREA). The algebras of this type can be associated to a
large class of braidings (solutions of the quantum Yang-Baxter
equation) of the Hecke type.

In the present paper we consider some quotients of the mREA --- the
non-commutative (NC) orbits. The terminology is motivated by a close
connection of these quotients with the coordinate algebras of orbits
in $gl(n)^*$.

For each generic NC orbit $\A$ we establish certain combinatorial
relations among its central elements. The relations can be
interpreted as the higher NC counterparts of the Newton identities.
Recall, that the classical Newton relations connect the elementary
symmetric functions of some commutative variables and the power
sums of the same variables. Besides, we construct a large family of
projective $\A$-modules (in what follows all modules are assumed to
be finitely generated and one-sided).

Also, we introduce an algebra $\KK(\A)$ which is an analog of
$K^0(\Fl(\C^n))$, where $\Fl(\C^n)$ is a flag variety, and define a
$q$-analog of the Euler characteristic of line bundles on the flag
variety which is is well defined on $\KK(\A)$.

A particular case of the NC orbits is the set of the so-called
quantum orbits arising from the quantization  of a certain Poisson
pencil on semisimple orbits in $gl(n)^*$ (i.e., $GL(n)$ orbits of
semisimple elements of $gl(n)^*$). The Poisson pencil is generated
by the two brackets: the Kirillov bracket and another one, related
to a classical $r$-matrix. One of the main properties of the algebra
$\A$ arising from  the quantization of the Poisson pencil is that
the product $\mu: \A\otimes \A \rightarrow \A$ is equivariant
(covariant) with respect to the action of the quantum group  $\Uq$
\be
U\triangleright\mu(a\otimes
b)=\mu(U_1\triangleright a\otimes U_2\triangleright b),\quad a,b
\in \A,\quad U\in\Uq, \label{equi}
\ee
where $\De(U)=U_1\ot U_2$. In this case (the $\Uq$ case for short)
at the limit $q\rightarrow 1$ we get the $SL(n)$-equivariant (or,
if we pass to the compact form, $SU(n)$-equivariant) algebras which
are also called {\em the fuzzy orbits}.

Let us describe the algebras in question in more detail. We start
from the definition of the reflection equation algebra connected
with a braiding of the Hecke type. The initial data for its
construction is a braiding $R$
$$
R:\quad\vv\to\vv,
$$
which is a solution of the quantum Yang-Baxter equation (\ref{YBE})
satisfying additionally the second order equation (\ref{Hecke}). Here
V is a finite dimensional vector space, $\dim V = n$. Such a braiding
will be called {\em the Hecke symmetry}. A well-known example is
connected with the quantum group $\Uq$ when the corresponding Hecke
symmetry is the image of the universal $R$-matrix in the fundamental
vector representation of $\Uq$.

In \cite{G} there were constructed other examples of Hecke symmetries
such that the Hilbert-Poincar\'e series of the associated "symmetric"
and "skew-symmetric" subalgebras of the tensor algebra $T(V)$ differ
from the classical ones. Below we shall additionally assume the
Hilbert-Poincar\'e series of the "skew-symmetric" subalgebra to be a
monic polynomial. Such a Hecke symmetry will be called {\em even} and
the degree $p$ of this polynomial will be called {\it the symmetry
rank} of $R$.

Consider now an associative unital algebra generated by $n^2$
indeterminates $l_i^{\,j}$, $1\le i,j\le n$ satisfying the following
relation
\be
R\, L_1\,R \,L_1 -  L_1\,R\, L_1\,R = \hbar(R\,L_1 - L_1\,R),
\qquad  L_1\equiv  L\otimes I,   \label{defREA}
\ee
where $\hbar$ is a formal parameter, $I$ is an $n\times n$ unit
matrix and $L =\|l_i^{\,j}\|$ is a matrix with entries $l_i^{\,j}$.
If $\h=0$, we call this algebra the {\it reflection equation
algebra} (REA) and denote it $\lq$; if $\h\not=0$, we call it
the {\it modified reflection equation algebra} (mREA) and denote
it $\lhq$.

In the $\Uq$ case, being specialized at $q=1$, the mREA coincides
with the enveloping algebras $U(gl(n)_\h)$ (hereafter, given a Lie
algebra $\gg$ with the bracket $[\,\cdot\,,\cdot\,]$, the symbol
$\gh$ will denote the Lie algebra with the bracket
$\h[\,\cdot\,,\cdot\,]$). In some sense the mREA is a "braided"
analog of the enveloping algebra $U(gl(n)_\h)$. For a motivation
of this treatment, see section \ref{sec:rea-br}.

Now we list some properties of the mREA. First of all, the
definition implies that the generators $l_i^j$ of any mREA obey the
quadratic-linear commutation relations. Second, the category $\lmod$
of {\em equivariant} finite dimensional representations of the mREA
corresponding to an even Hecke symmetry is close to the category
$U(gl(p))$--Mod, where $p$ is the symmetry rank of the Hecke
symmetry $R$ (see section \ref{sec:rea-br}). For example, simple
objects of the category $\lmod$ can be labelled by signatures
(partitions) $\la=(\la_1,\dots,\la_p)$, $ \la_1\geq\dots \geq
\la_p$, similarly to the category of $gl(p)$-modules. Besides, the
Grothendieck rings of these categories are isomorphic \cite{GLS1}.

To explain the term "equivariant" we restrict ourselves to the $\Uq$
case. Consider an  $\Uq$ module $V$. Then, as is well known, the
space $\End(V)$ of the internal endomorphisms of $V$ is endowed with
an $\Uq$-action, too. The mREA can as well be equipped with an
$\Uq$-action satisfying (\ref{equi}). In this case a representation
$\pi:\lhq\to\End(V)$ is called {\it equivariant} if $\pi$ commutes
(as a mapping) with the $\Uq$-action.

The next important property of the mREA is the existence of a monic
polynomial ${\cal CH}(x)$ of degree $p$ such that the matrix $L$
entering formula (\ref{defREA}) satisfies the Cayley-Hamilton (CH)
identity
$$
{\cal CH}(L) = \sum_{k=0}^p(-L)^{p-k}\sigma_k(L) \equiv 0,\qquad
\sigma_0(L)={\rm id}_{\cal L}
$$
where the coefficients $\sigma_k(L), \,0\le k\le p$ are linear independent
generators of the center $Z(\lhq)$
 of the
mREA. This identity was proved in \cite{GPS} for the non-modified
REA (the $\Uq$ case was previously considered in \cite{PS}). The
above CH identity for the mREA can be established by a linear change
of generators \cite{GS1}.

Let $\chi:\zl\to\K$ be a character of the center. The character
$\chi$ is completely fixed by its values on $\sigma_k$
$$
\chi(\sigma_k)=\sum_{1\leq i_1<\dots<i_k\leq p}\mu_{i_1}\dots
\mu_{i_k}, \quad 1\le k\le p
$$
where the numbers $\mu_i\in\K$ are assumed to be distinct. The
quotient of the algebra $\lhq$ modulo the ideal generated by the
elements $z-\chi(z)$, where $z\in\zl$, will be called an {\it NC
orbit} and denoted $\lhqc$ (although as explained in
remark \ref{DM} such an "NC orbit" in the $U_q(sl(n))$ case
can arise from a quantization of a union of some orbits in
$gl(n)^*$).

Observe that, being switched to the
algebra $\lhqc$, the CH identity for $L$ takes the form
$$
\prod_i(L-\mu_i \Id)=0.
$$
(Thus, the numbers $\mu_i$ are thought of as eigenvalues of the
matrix $L$ with entries from $\lhqc$.)

A great importance of the CH identity is occasioned by the fact that
it allows us to construct a family of idempotents from
${\rm Mat}(\lhqc)$ and therefore projective modules. Moreover, the
family of idempotents (and the corresponding projective modules) can
be essentially enlarged. For this purpose, we construct a series of
higher order matrices $\lm$ and polynomials
$\ch_{(m)},\,\,m=2,3,\dots$ with central coefficients such that the
higher order CH identities $\ch_{(m)}(\lm)=0$ are satisfied. We also
set $L_{(1)}=L$, $\ch_{(1)}=\ch$.

Upon restricting to the NC orbit $\lhqc$, we come to the series of
polynomials ${\cal CH}_{(m)}^{\chi}$ with the numerical
coefficients. Assuming the roots of these polynomials to be distinct
for any $m$ we find $p_m = \deg\,\ch_{(m)}$ idempotents $e_\k(m)\in
{\rm Mat}(\lhqc)$, each corresponding to a projective
$\lhqc$-module.

One of the main aims of this paper is to prove the existence of the
higher order CH identities and to compute the coefficients of the
polynomials $\ch_{(m)}$. This is rigorously done for the mREA
associated with any even Hecke symmetries of rank $p=2$. For
$p\geq 3$ we present an explicit formula for these polynomials as a
conjecture.

Our other aim is to compute the values of the central elements
$\Tr_R\lm^s$, $s\ge 1$, for a generic NC orbit. Here $\Tr_R:
{\rm Mat}(\lhq)\rightarrow \lhq$ is the trace defined by the initial
Hecke symmetry $R$. It is closely related to the {\it categorical
trace} which is discussed in section \ref{sec:rea-br} (see
\cite{GLS1} for more detail). For example, in the $\Uq$ case $\Tr_R$
is the well known quantum trace which is a weighed sum of the
diagonal entries, but in the general case it can be more
complicated. We express $\Tr_R \lm^s$ (up to the aforementioned
conjecture) in terms of the eigenvalues of the matrix $L$. For $m=1$
we treat these expressions as a parametric resolution of the Newton
relations. For $m>1$ they are thought of as higher analogs of the
Newton relations.

As a byproduct we compute the value of $\Tr_R \,e_\k(m)$ on generic
NC orbits. In contrast with the usual orbits when these quantities
are always equal to 1 (since the related projective modules
correspond to line bundles), for the NC orbits they are more
informative. We show that the assignment $M_\k(m) \mapsto \Tr_R
e_\k(m)$, where the module $M_\k(m)$ corresponds to the idempotent
$e_\k(m)$, can be considered as an analog of the Euler
characteristic of a line bundle over a flag variety. Developing the
analogy with a flag variety, we can introduce a multiplicative
structure on the set of the modules $M_\k(m)$ and construct an
algebra $\KK(\lhqc)$ playing the role of $K^0$ of the flag variety.

Completing the Introduction, we consider the $\Uq$ case in more
detail. As we have already mentioned, in this case the NC orbits
arise from a quantization of a Poisson pencil on some algebraic
varieties. So, it is natural to consider the problem of quantization
of the vector bundles over these varieties in terms of projective
modules. In this case the projective modules $M_\k(m)$ over the
quantum orbits are nothing but deformations of the corresponding
modules over the coordinate rings of the initial varieties.

The problem of quantization of semisimple orbits in $\gg^*$ was
considered in numerous papers. In \cite{DM1} the full solution of
this problem was given in terms of the mREA and its appropriate
quotients related to the minimal polynomials of orbits in question.
In \cite{DM2} this method of quantization is compared in the
$U(sl(n))$ case with an approach based on the so-called generalized
Verma modules. In particular, the authors produce a formula
describing the eigenvalues $\mu_i$ to be functions of the
corresponding generalized Verma module weight 
$\lambda_i$\footnote{In fact, this formula is well-known, cf. for
example, \cite{BR} where the eigenvalues of the Casimir elements
$\Tr L^s$ in finite dimensional representations are computed.}.
We find a q-analog of this formula by using different methods. (Note
that in the $U(sl(n))$ case this formula is obtained by means of the
coproduct, see section 2, whereas in the mREA we did not find any
convenient coproduct.)

The paper is organized as follows. In the next section we consider
the NC orbits arising from the quantization of the Kirillov bracket
on some semisimple orbits in $gl(n)$ (or $su(n)$). We call them {\em
the fuzzy orbits}, the term is motivated by numerous papers devoted
to "the fuzzy physics". In contrast with all those papers, we
present a general scheme of constructing a large family of
projective modules via the CH identities. In subsequent sections we
generalize this scheme to the mREA associated with a large class of
Hecke symmetries and thereby we construct a similar family of
projective modules for their appropriate quotients
("$q$-fuzzy orbits").

Note that the methods of section \ref{sec:fuzzy} are close to those
going back to the pioneering paper \cite{Ko} and having been used in
numerous works devoted to the {\it characteristic identities} (which
are nothing but a specialization of the CH identities to concrete
representations of the algebras in question, cf. \cite{Go}, \cite{O}
and references therein\footnote{Note that there are known $q$-analogs
of such characteristic identities related to the quantum groups (cf.
\cite{GZB}). However, the quantum groups are not deformations of
commutative algebras and are not convenient objects for constructing
projective modules and the related Newton identities. A set of
specific New\-ton-Cay\-ley-Ha\-mil\-ton identities exists in
algebras dual to the quantum groups, cf. \cite{IOP}. However, they
are useless for constructing the projective modules over these
algebras.}).

In section \ref{sec:rea-def}, we introduce the REA and list its
basic properties in detail.

In section \ref{sec:rea-br}, we give reasons allowing us to treat
the mREA as a braided analog of the enveloping algebra.

Section \ref{sec:5} is devoted to the proof of a parametric
resolution of the {\em basic Newton identity}.

In section \ref{sec:6} we give a proof of the relation between
eigenvalues of the matrix $L$ and those of the highest matrices
$L_{(m)}$ in the particular case $p=2$. At the end of the section
we compute the values of the central elements $\Tr_R \lm^k$. As a
byproduct we get the aforementioned formula featuring the 
eigenvalues of the matrix $L$ as a function of a partition $\la$.

Finally, in section \ref{sec:q-eiler} we introduce the algebra
$\KK(\lhqc)$ and define a $q$-analog of the Euler characteristic.

{\bf Acknowledgement} The authors are grateful to the
Max-Planck-Institut f\"ur Mathematik (Bonn) where the paper was
started and the Institute Mittag-Leffler (Stockholm) where the paper
was completed for the warm hospitality and the stimulating
atmosphere. One of the authors (D.G.) is grateful to A. Mudrov for
valuable discussions.

\section{Fuzzy orbits}
\label{sec:fuzzy}

In this section we are dealing with algebras arising from the
quantization of the Kirillov bracket alone. First, we describe our
initial object --- the coordinate algebra of a generic classical
orbit. Let us fix a diagonal matrix
\be
M={\rm diag}(\mu_1,\dots , \mu_n),\quad \mu_i\in \C
\label{elem}
\ee
with simple $\mu_i$. We treat this matrix as an element of
$\gg^* = gl(n)^*$ (identified with $\gg=gl(n)$) and consider its
orbit $\Om$ with respect to the $GL(n)$-action
\be
GL(n)\ni g:\quad M \mapsto M_g=g^{-1}\, M\, g.
\label{act}
\ee
The orbit $\Om$ is a closed affine algebraic variety and its
coordinate algebra $\K(\Om)$ can be described as follows. Let
$\K(\gg^*)$ be the polynomial algebra in $n^2$ commutative
indeterminates $l_i^j,\,1\le i,j \le n$, which are nothing but the
coordinate functions. Then we have $\K(\Om)=\K(\gg^*)/{\cal I}$
where ${\cal I}\subset \K(\gg^*)$ is the ideal generated by the
elements
\be
\Tr\,L^k-\beta_k,\quad\,k=1,2, \dots , n.
\label{lk}
\ee
Here the matrix $L=\|l_i^j\|$ is composed of the indeterminates
$l_i^j$ ($i$ labels the rows and $j$ --- the columns) and
\be
\beta_k=\sum_{i=1}^n \mu_i^k,\quad k=1,2, \dots , n.
\label{alp}
\ee

\begin{remark}
\label{rem:nongener}
{\rm
Note that if $\Om$ is semisimple but not generic (i.e., the
eigenvalues $\mu_i$ of a given diagonal matrix $M$ are not simple)
the coordinate algebra is the quotient $\K(\Om)=\K(\gg^*)/
{\cal I}'$, where the ideal ${\cal I}'$ is of the following form.
Let
$$
{\cal P}(M) = (M-\mu_1I)\dots (M-\mu_rI)
$$
be the minimal polynomial of the matrix $M$ ($I$ stands for the
$n\times n$ unit matrix). Then $n^2$ entries of the matrix
${\cal P}(L)$ are polynomials in generators $l_i^j$. The ideal
${\cal I}'$ is generated by these polynomials and by elements
(\ref{lk}) with $k=1,\dots ,r-1$. Note, that if we disregard the
latter elements we get a union of all semisimple orbits possessing
the same minimal polynomial.}
\end{remark}

If the initial matrix $M$ is such that $\Tr\, M=0$, then the
corresponding orbit is embedded into $sl(n)^*$. Since
$$
\K(sl(n)^*)=\K(gl(n)^*)/\{\Tr\, L\},
$$
the coordinate algebra of the corresponding orbit can be realized as
above but with $\beta_1=0$. (Hereafter $\{S\}$ stands for the ideal
generated by a set $S$.)

If all eigenvalues of the matrix $M$ are real, we can consider the
matrix ${\ii} M$ (here ${\ii}=\sqrt{-1}$) as an element of $u(n)^*$
(or $su(n)^*$ if $\Tr\, M=0$). Choosing the generators $x_i^j$ for
$i\le j$ and $y_i^j$ for $i<j$ such that
$$
l_i^j=x_i^j+ {\ii} y_i^j\quad{\rm for}\;\;i<j, \qquad
l_i^j=x_j^i-{\ii} y_j^i\quad{\rm for}\;\;i>j, \qquad{\rm and}
\qquad
l_i^i=x_i^i,
$$
we get a compact real variety which is an $SU(n)$-orbit (i.e., in
(\ref{act}) we assume that $g\in SU(n)$). Consequently, we consider
its coordinate algebra as an $\R$-algebra.

Now let us pass to the fuzzy orbits. In the spirit of "NC affine
algebraic geometry" we realize them via some explicit relations on
generators. Consider again the matrix $L=\|l_i^j\|$, but now we let
the generators $l_i^j$ to satisfy the defining relations of the
algebra $U(\gh)$ with $\gg=gl(n)$:
$$
l_i^j\,l_m^n-l_m^n\,l_i^j-\h(l_i^n\,\de_m^j-l_m^j\,\de_i^n)=0.
$$
Then the matrix $L\in\Mat_n(U(\gh))$ obeys a polynomial relation
\be
\ch(L)=\sum_{k=0}^n(- L)^{n-k}\sigma_k( L) = 0,
\label{CH1}
\ee
such that the coefficients $\sigma_k( L)$ are central and
$\sigma_0( L) = 1$. This fact is well known. An expression of the
coefficients $\sigma_k( L)$ in terms of the generators $l_i^{\,j}$
can be found in \cite{GS1}. Below we present them in a convenient
form for a more general case of NC orbit (see section \ref{sec:5}).

It is also well known that the center $Z(U(\gh))$ of the algebra
$U(\gh)$ is generated by $\sigma_k( L)$ for $1\le k\le n$. Another
family generating the center is $s_k(L)=\Tr L^k$, $ 1\le k \le n$.
Therefore, any character
$$
\chi:\; Z(U(\gh))\to \K
$$
is completely determined by its values on the generators of the
center $\chi(\sigma_k( L)) = \al_k$ or $\chi(s_k(L))=\beta_k$, $k=
1,\dots,n$.

Consider the quotient algebra
\be
\lhc = U(\gh)/\{z-\chi(z)\;|\; z\in Z(U(\gh)) \},
\label{alg}
\ee
where $\chi$ is a fixed character. Being  switched  to $\lhc$,
relation (\ref{CH1}) takes the form
\be
\chc(L)\stackrel{\mbox{\tiny def}}{=}\sum_{k=0}^n(- L)^{n-k}
\chi(\sigma_k(
L)) = \sum_{k=0}^n(- L)^{n-k}\al_k=0.
\label{CH2}
\ee
In what follows the superscript $\chi$ means that we have passed
from $U(\gg_\h)$ to $\lhc$ (and similarly, for other algebras
below).

Also, consider the associated numerical equation
\be
\sum_{k=0}^n(- \mu)^{n-k}\al_k = 0. \label{CH3}
\ee
Denoting the roots of this equation by $\mu_i$ we get
\be
\alpha_k = \chi(\sigma_k) =
\sum_{1\le i_1<i_2<\dots <i_k\le n}\mu_{i_1}\mu_{i_2}\dots
\mu_{i_k}.
\label{gen}
\ee
Hereafter (up to the last section) we assume the numbers
$\mu_i\in\K$ to be fixed and the character $\chi$ to be defined
by the set of values (\ref{gen}).

\begin{definition}
\label{def:fuz-orb}
The roots of equation (\ref{CH3}) will be called the eigenvalues of
the matrix $L$ on the orbit $\lhc$ (that is, when entries of $L$
belong to the algebra $\lhc$). The  algebra $\lhc$ will be called
the 1-generic fuzzy orbit if the eigenvalues $\mu_i$ are simple.
\end{definition}

In what follows we shall only consider the 1-generic fuzzy orbits
without specifying it each time. Note, that the 1-generic fuzzy
orbits can arise from a quantization of semisimple but not generic
orbits (see remark \ref{DM} below).

For any fuzzy orbit, we introduce $n$ idempotents:
\be
e_j=\prod_{i\not=j}{{(L-\mu_i\Id)}\over{(\mu_j-\mu_i)}},\qquad
j=1,\dots,n.  \label{proj1}
\ee
Identity (\ref{CH2}) leads to
$$
e_i\,e_j=\delta_{ij}\,e_i,\qquad \sum_{i=1}^n e_i=\Id.
$$
The CH identity (\ref{CH1}) and all related objects (idempotents
(\ref{proj1}), corresponding projective $\lhc$-modules, etc.) will
be called {\em basic}.

Now, we shall describe a regular way of constructing some
{\em higher} analogs of these objects. To this end, consider the
category of finite dimensional representations of the algebra
$U(\gg)$. Its simple objects $\vl$ are labelled by sequences of
numbers
\be
\la=(\la_1,\dots ,\la_{n}),
\label{signature}
\ee
where $\la_i-\la_{i+1}$ are nonnegative integers. Following \cite{W}
we call these sequences the {\em signatures}. If moreover, $\la_{i}$
themselves are nonnegative integers and $\sum \la_i=m$ we call the
signature $\lambda$ {\it ordered partition} of the integer $m$.
Since the algebras $U(\gh)$ and $U(\gg)$ are isomorphic to each
other, the objects $\vl$ can be equipped with an $U(\gh)$-action.

Consider a left $U(\gh)$-module $\vl$ and let
$$
\pi_{\la}:\; U(\gh)\to\End(\vl)
$$
be the corresponding left irreducible representation. All
representations in question are assumed to be {\em equivariant},
i.e., they commute with the $GL(n)$-action, where we suppose
$U(\gh)$ to be equipped with the adjoint $GL(n)$-action.

Introduce now the map
$$
\pi_\la^{(2)}=\Id\ot\pi_\la :\quad U(\gh)\ot U(\gg)\to U(\gh)\ot
\End(\vl)
$$
(note, that in the second factor we put $\h=1$). On fixing a basis
in the space $\vl$, we can identify the spaces $\End(\vl)$ and
$\Mat_{\nl}(\K)$, where $\nl=\dim\vl$. Consequently, the spaces
\be
U(\gh)\ot \End(\vl)\quad
{\rm and}\quad U(\gh)\ot \Mat_{\nl}(\K)=\Mat_{\nl}(U(\gh))
\label{ident}
\ee
can be identified. Therefore, the above map $\pi_\la^{(2)}$ is of
the form
$$
\pi_\la^{(2)}:\quad U(\gh)\ot U(\gg)\to \Mat_{\nl}(U(\gh)).
$$
\begin{remark}
\label{rem:3}
{\rm
Note that the spaces $V_\lambda$ and $V_\mu$ whose signatures are
different by a constant shift
$$
\lambda_i -\mu_i = z,\quad 1\le \forall\, i\le n,\quad
z\in{\Bbb K}
$$
have equal dimensions $n_\lambda  = n_\mu$. Making the
transformation
\be
\pi_\lambda(e_i^{\,j})\mapsto \pi_\lambda(e_i^{\,j}) -
z\,\delta_i^{\,j}
{\rm id}_{V_\lambda}\quad \forall\,z\in{\Bbb K}
\label{ed-sdv}
\ee
where $e_i^{\,j}$ are the $U(gl(n))$ generators, one can convert the
$U(gl(n))$-representation $\pi_\lambda$ into $\pi_\mu$. In
particular, setting $z=\lambda_n$ we come to the representation
$\pi_{\hat
\lambda}$ where $\hat\lambda$ is the {\it partition} of the form
\be
\hat \lambda = (\hat \lambda_1, \dots, \hat \lambda_{n-1},0),
\quad \hat\lambda_i = \lambda_i-\lambda_n.
\label{lambda-hat}
\ee
Note that the objects $V_\lambda $ and $V_{\hat \lambda}$ differ as
$U(gl(n))$-modules but coincide as $U(sl(n))$-ones (recall that
simple $U(sl(n))$-modules are labelled by partitions 
(\ref{lambda-hat})). Therefore,
$\dim V_\lambda = \dim V_{\hat \lambda}$.

Moreover, to the space $V^*_\lambda$ (dual to $V_\lambda$) there
corresponds the space $V_{\lambda^{\!*}}$ in the category of finite
dimensional $U(gl(n))$-modules. Here the signature $\lambda^{\!*}$
is defined as
\be
\lambda^{\!*} = (-\lambda_n,\dots,-\lambda_1).
\label{lambda-star}
\ee
Note that in the category of $U(sl(n))$-modules the dual object is
labelled by the partition $\hat{\la^{\!*}}$. }
\end{remark}

Now we apply $\pi_\la^{(2)}$ to the {\it split Casimir element}
\be
\Cas= l_i^j\ot l_j^i.
\label{cas-kl}
\ee
Hereafter the summation over the repeated indices is understood. For
$\h= 1$, the image of the split Casimir element  under the product
map $U(\gg)^{\ot 2}\to U(\gg)$ becomes the usual quadratic Casimir
element $s_2=\Tr L^2\in U(\gg)$.

Denote  the matrix transposed to $\pi^{(2)}_\la(\Cas)$ by
$L_{\la}$, i.e.
$$
L_{\la}^t = \pi^{(2)}_\la(\Cas) = l_i^j\ot \pi_\la(l_j^i).
$$
It can be easily seen that if $\la_0=(1,0,\dots ,0)$, then
$L_{\la_0}=L$ provided that the basis $\{x_i\}\in V$ is fixed in
such a way that $\pi_{\la_0}(l_i^j)\triangleright x_k =
\delta_k^j x_i$.

We emphasize,  that the element $\pi_{\la}^{(2)}(\Cas)$ actually
belongs to the algebra $\Mat^{\Inv}_{\nl}(U(\gh))$  --- the
$GL(n)$-invariant subalgebra of $\Mat_{\nl}(U(\gh))$. This follows
from the fact that we are working with equivariant representations
$\pi_\lambda: U(\gh) \to \End(\vl) = \Mat_{\nl}(\K)$. Here we assume
that the action of $GL(n)$ on $\Mat_{\nl}(U(\gh))$ is defined on the
base of identification (\ref{ident}). In \cite{K} the algebras
$\Mat^{\Inv}_{\nl}(U(\gh))$ were called {\it the family algebras}
(classical if $\h=0$ and quantum if $\h\not=0$), see also \cite{R}.

\noindent
{\bf Example}\ \  Let $n=2$. Then
$$
\Cas=a\ot a+b\ot c +c \ot b+ d\ot d,\quad {\rm where} \quad a=
l_1^1, \,\,
b=l_1^2,\,\, c=l_2^1,\,\,d=l_2^2.
$$
Set $\la=(l+1,l),\,\,l\in\K$ and consider the representation
$\pi_\la$ of $U(gl(2))$. We have
$$
\pi_\la(a)=\left(\matrix{1+l&0\cr
0&l}\right),\;\pi_\la(b)=\left(\matrix{0&1\cr 0&0}\right),\;
\pi_\la(c)=\left(\matrix{0&0\cr
1&0}\right),\;\pi_\la(d)=\left(\matrix{l&0\cr 0&1+l}\right).
$$
The corresponding matrix $L_\la$ is
$$
L_\la=\left(\matrix{a&b\cr c&d}\right)
+l(a+d)
\left(\matrix{1&0\cr 0&1}\right).
$$
If $l=0$, we have just the matrix $L$.

Below we shall restrict ourselves to the matrices $L_\la$
corresponding to the partitions $\la=(m)=(m,0,0,\dots ,0)$. These
matrices, as well as related objects, will be denoted $\lm$,
$\pi_{(m)}$, $\vm$, etc.

The matrix $\lm$ as well as $L=L_{(1)}$ satisfies the CH identity.
Now we describe the construction of the corresponding CH polynomial.

In what follows we shall often use the set of all possible
(non-ordered) partitions $\k$ of the integer $m$ with the length not
greater than $n$. By definition, $\k$ is the set of $n$ integers
$k_i$ with the following properties
\be
\k=(k_1,\dots ,k_n),\quad k_i\geq 0,\quad |\k|=k_1+\dots +k_n=m.
\label{part}
\ee
For any partition $\k\vdash m$ from the set (\ref{part}), let
\be
\mu_\k(m)=\sum_{i=1}^n k_i\mu_i+\h\sum_{1\le i<j\le n} k_i k_j
\label{high}
\ee
where $\mu_i$, $1\le i\le n$, are elements of the algebraic closure
of the center $Z(U(\gg_\h))$. They are the solutions of the set of
polynomial relations
\be
\sum_{1\le i_1<i_2<\dots <i_k\le
n}\mu_{i_1}\mu_{i_2}\dots \mu_{i_k}=\sigma_k(L),
\quad 1\le k\le n,
\label{elsym}
\ee
where $\sigma_k(L)$ are the coefficients of the basic CH polynomial
(\ref{CH1}).

Let us define the monic polynomial
\be
\ch_{(m)}(x) = \prod_{k\vdash m}(x-\mu_\k(m)).
\label{pol-form}
\ee
Since its coefficients are symmetric functions in $\mu_i$ they can
be expressed via the elementary symmetric functions. By virtue of
(\ref{elsym}) we conclude that the coefficients of $\ch_{(m)}$
belongs to $Z(U(\gg_\h))$. Besides,
$\deg \ch_{(m)}(x) = n_m=\dim\vm$.

\begin{proposition}
\label{prop:3}
The polynomial $\ch_{(m)}(x)$ is a CH polynomial for $\lm$, namely,
we have
$$
\ch_{(m)}(\lm)\equiv 0.
$$
\end{proposition}

When switching to the 1-generic orbit $\lhc$ (\ref{alg}), we fix the
eigenvalues $\mu_i$ of the matrix $L$ (see (\ref{gen})) and thereby
fix the values of the quantities $\mu_\k(m)$. Similarly to the case
$m=1$, the quantities $\mu_\k(m)$ are called {\it eigenvalues} of
the matrix $\lm$ {\it on the orbit $\lhc$}. Besides, for the sake of
uniformity, we put $\mu_\k(1)=\mu_i$ for $|\k|=1$ and
$k_j= \delta_{i,j}$.

\medskip \noindent
{\bf Proof\ \ } In what follows, besides the representations
$\pi_\la$, we also need the right representations $\opi_\la$. The
representation $\opi_\la$ is defined in the space $V_\la^{\!*}$ dual
to $V_\la$
$$
\langle X\triangleleft\opi_\la(a), Y\rangle=\langle X,\pi_\la(a)
\triangleright Y  \rangle,\quad X\in V_\la^{\!*},\; Y\in V_\la,\;
a\in U(\gh).
$$
Hereafter, $\triangleleft$ (resp., $\triangleright$) stands for the
right (resp., left) action of a given operator on a given vector.
Replacing the operators $\opi_\la(a)$ by $-\opi_\la(a)^*$ where $^*$
is an involution, we get a left representation $\pi_{\la^{\!*}}$ of
the algebra $U(\gh)$ which is called contragradient to $\pi_\la$ and
is labelled by the signature $\lambda^{\!*}$ defined in
(\ref{lambda-star}) (cf. \cite{W}).

The idea of our proof consists in the following. To verify a
relation in the algebra $U(\gh)$, we consider the image of this
relations in an arbitrary representation $\opi_\la$ and prove that
this relation is true in all such representations. Then, since any
element of the enveloping algebra with trivial image in any
representation  $\opi_\la$ is trivial (cf. \cite{D}), we conclude
that our relation is valid in the algebra itself. Moreover, it
suffices to check the relation in question for $\h= 1$, since by
rescaling we can pass to an arbitrary $\h\not=0$.

Let us introduce the following notation:
$$
\llm=\opi_\la(\Lm),
$$
where applying a representation $\opi_\la$ to a matrix means
applying it to any entry of  the matrix in question. Thus,
$\opi_\la(\lm)$ is a matrix with entries from $\End(\vl)$ which in
turn is identified with $\Mat_{n_\la}(\K)$.

So, given a fixed $m$ and an arbitrary $\lambda$, we should prove
the identity $\ch_{(m)}(\llm)\equiv 0$ with the polynomial
$\ch_{(m)}(x)$ defined in (\ref{pol-form}).

For a fixed $m$, let us first consider a subset of representations
$\opi_\la$ corresponding to signatures
$\la=(\la_1,\la_2,\dots ,\la_n)$ such that
\be
\la_i- \la_{i+1}\geq m\quad {\rm for}\quad 1\le i\le n-1.
\label{proper}
\ee

For such signatures,  irreducible components in the tensor product
$V_\la^*\ot\vm$ are of the form $V_{\la^{\!*}+\k}$, where $\k$ runs
over the full set of all possible partitions (\ref{part}). The sum
$\la+\k$ means the signature with the components
\be
(\la+\k)_i = \la_i+k_i.
\label{sum-lk}
\ee
Since all $k_i\le m$, then, due to restriction (\ref{proper}), we
have $(\la^{\!*}+k)_{i}\geq (\la^{\!*}+k)_{i+1}$ and the space
$V_{\la^{\!*}+\k}$ is well-defined.

Now, assuming $\la$ to satisfy (\ref{proper}), we use the split
Casimir element (\ref{cas-kl}) in order to calculate the eigenvalues
of the matrices $L_{(\la, 1)}$ and $\llm$. With this purpose we
insert the split Casimir $\Cas\in U(\gg)\ot U(\gg)$ between the
factors of the tensor product $V_\la\ot \vm$ and consider the image
\be
(\opi_\la\ot \pi_{(m)})(V_\la^*\ot\Cas\ot \vm)
\stackrel{\mbox{\tiny def}}
{=}
\sum_{1\le i,j\le n} (V_\la^* \triangleleft \opi_\la(l_i^j))
\ot (\pi_{(m)}(l_j^i) \triangleright\vm).
\label{cas-lr}
\ee
Thus, we have represented $\Cas$ as a linear operator $\Caslm$ in
the space $V_\la^*\ot \vm$
$$
\Caslm:\quad V_\la^*\ot \vm\to V_\la^*\ot \vm.
$$
This operator commutes with the action of the group $GL(n)$ and
therefore it is scalar on each irreducible component
$V_{\la^{\!*}+\k}\subset V_\la^*\ot \vm $. Since in an appropriate
basis of $V_\la^*\ot \vm$ the matrix of the operator $\Caslm$ 
coincides with $\llm$, then the eigenvalues of the matrix
$L_{(\la, m)}$ are nothing but the eigenvalues of the operator
$\Caslm$.

Let $\mu_\k(\la,m)$ be the eigenvalue of $\Caslm$ on the component
$V_{\la^{\!*}+\k}$ (if $m=|\k|=1,\,k_j=\delta_{i,j}$ we shall also
use the notation $\mu_i(\la)$ for the eigenvalues of
$\Cas_{(\la,1)} = L_{(\la,1)}$). One can prove that
$$
\mu_\k(\la,m)=-{{1}\over{2}}(s_2(\la^{\!*}+\k)-s_2(\la^{\!*})
-s_2((m))),
$$
where $s_2(\la)$ stands for the value of the quadratic Casimir
element on the irreducible module $\vl$. This formula is an
immediate consequence of the Leibnitz rule (or in other words, the
coproduct in the algebra $U(\gg)$). The negative sign appears due to
the passage from the representation $\opi_\la$ to $\pi_{\la^{\!*}}$.

The straightforward calculation gives
$$
s_2(\la)=\sum_{1\le i\le n} (\la_i^2+ \la_i(n+1-2i)).
$$

Therefore, for any  $\la$ satisfying (\ref{proper}), we have
\be
\mu_\k(\la, m)=-\sum_{i=1}^n k_i(-\la_{n-i+1}-i+1)+
\sum_{1\le i<j\le n} k_i k_j.
\label{muk}
\ee
In particular, if $m=1$, we get the eigenvalues of $L_{(\lambda,1)}$
\be
\mu_i(\la)=\la_{n-i+1}+i-1,\quad i=1,\dots ,n.
\label{mu}
\ee
To pass to a generic $\h$, we should multiply all these eigenvalues
by $\h$. Finally, we conclude that for  all $\h$ and the signatures
$\la$ satisfying (\ref{proper}) the eigenvalues of $\llm$ and
$L_{(\la,1)}$ are indeed connected by (\ref{high}):
\be
\mu_\k(\la, m)=\sum_i k_i \mu_i(\la)+ \h\,\sum_{1\le i<j\le n}
k_i k_j.
\label{mu-lam}
\ee

So, if we apply such a representation $\opi_\la$ to the matrix
$\ch_{(m)}(\lm)$ we obtain 0. Now we must get rid of restriction
(\ref{proper}).

For this purpose,  observe  the following. Let $m=1$. Then once
$\la_i\geq \la_{i+1}+1$ the eigenvalues $\mu_i(\la)$ are given by
(\ref{mu}). However, this formula is valid even if the condition
(\ref{proper}) with $m=1$ is not fulfilled. Indeed, formula
(\ref{reduc}) below gives correct values of central elements
$\Tr L^s_{(\la, 1)}$ if we assume the eigenvalues $\mu_i(\la)$ to be
always given by (\ref{mu}) without any restrictions on $\la$. It can
be verified through calculating $\Tr L^s_{(\la, m)}$ by other means.
For instance, it can be done as in \cite{BR}\footnote{\label{foot:2}
In order to compare the formula from \cite{BR} with our result it
suffices to invert the numbering of the eigenvalues $\mu_i$.} or
$\Tr L^s_{(\la, 1)}$ can be obtained from the $q$-Newton identities
below. By inverting the Newton identities we can express the central
elements $\sigma_k(L)$ in terms of $\Tr L^s$. This implies that
formulae (\ref{mu}) are true  for any representation $\opi_\la$.

Now, we are able to complete the proof. Even if $\la$ does not
satisfy (\ref{proper}) for a given $m$, the eigenvalues of the split
Casimir $\Caslm$ form a subset of the set $\{\mu_\k(m)\}$ with
$\mu_\k(m)$ defined in (\ref{high}). So, if we apply the
representation $\opi_\la$ to the matrix $\ch(L_{(m)})$ we get 0.
This implies the statement. \hfill\rule{6.5pt}{6.5pt}

\begin{remark} \label{DMO}
{\rm
As we have seen above, the degree of the polynomial $\ch_{(\la,1)}$
is always equal to $n$, even though the number of the eigenvalues of
the split Casimir $\Cas_{(\la,1)}$ can be smaller. This means that
not all eigenvalues $\mu_i(\la)$ of the matrix $L_{(\la,1)}$ can be
found from the split Casimir but only those appearing in the minimal
polynomial (see remark \ref{DM}). However, formula (\ref{reduc})
below is always true since extra eigenvalues which do not come in
the split Casimir spectrum have vanishing coefficients in the r.h.s.
of (\ref{reduc}).}
\end{remark}

\begin{definition}
For a given natural $m$, a fuzzy orbit will be called $m$-generic if
it is 1-generic and the eigenvalues $\mu_\k(m)$ are simple. A fuzzy
orbit will be called generic if it is $m$-generic for any $m$.
\end{definition}

For any $m$-generic fuzzy orbit, we define
$n_m=\dim \vm=\left(\matrix{n+m\cr m}\right)$ idempotents in the
usual way:
\be
e_\k(m)=\prod_{\k'\not=\k}{\frac{(L_{(m)}-\mu_{\k'}(m)\Id)}
{(\mu_\k(m)-
\mu_{\k'} (m))}}.
\label{high-idemp}
\ee

The following proposition is an easy corollary of the CH identity
$\chc_{(m)}(\lm)=0$ and it is treated as an NC version of the
spectral decomposition.

\begin{proposition} \label{mult}
For any $m$-generic fuzzy orbit $\lhc$ characterized by eigenvalues
$\mu_i,\,\,1\le i\le n,$ there exist quantities $d_\k(m)$ such that
$$
\Tr\Lm^s=\sum_{|\k|=m}\mu_\k(m)^s\, d_\k(m), \quad s=1,2,\dots,
$$
where $\mu_\k(m)$ are given by (\ref{high}). Moreover, we have
$$
d_\k(m)=\Tr e_\k(m).
$$
\end{proposition}

The quantities $d_\k(m)$ will be called {\it the quantum
multiplicities}. Their precise values on fuzzy orbits are given by
the following proposition.

\begin{proposition}
\label{prop:6}
For any $m$-generic orbit $\lhc$ characterized by eigenvalues
$\mu_i,\,\,1\le i\le n$ we have
\be
d_\k(m)=\prod_{1\le i<j\le n}{{\mu_i-\mu_j-(k_i-k_j)\h} \over
{\mu_i-\mu_j}}. \label{razmer}
\ee
Equivalently, the following "higher  Newton identities" are valid
\be
\Tr \lm^s=\sum_{|\k|=m}(\sum_{i=1}^n k_i\mu_i+\h\sum_{1\le i<j\le
n}
k_ik_j)^s \prod_{1\le i<j\le n}{{\mu_i-\mu_j-(k_i-k_j)\h}
\over{\mu_i-\mu_j}}.
\label{princip}
\ee
\end{proposition}

\medskip\noindent
{\bf Proof\ \ } The proof is based on the following observation. It
is not difficult to see that the quantities $\Tr L_{(m)}^s$ in
(\ref{princip}) are (symmetric) polynomials in $\mu_i$. To define a
polynomial unambiguously, it suffices to fix it at a finite set of
values of its arguments.

As such a set we take the values $\mu_i(\lambda)$ (\ref{mu})
connected with the finite dimensional representations
parametrized by signatures $\lambda$. However, a subtle point here
is that not any representation $\opi_\lambda$ is admissible for our
purpose. In the proof we use the split Casimir element, therefore we
should take a signature $\lambda$ in such a way that the
corresponding eigenvalues $\mu_\k(\lambda,m)$ (\ref{mu-lam}) would
be all present in the spectrum of $\Caslm$ and, besides, they would
be all distinct. A simple analysis of the structure of
$\mu_i(\lambda)$ and $\mu_\k(\lambda,m)$ shows that these conditions
can always be met if the differences $\lambda_i-\lambda_{i+1}$ are
sufficiently large. Let us introduce the corresponding notion.

\begin{definition}
\label{def:m-ad}
A representation $\opi_\la$ (and  the corresponding signature $\la$)
will be called $m$-ad\-mis\-sible if condition (\ref{proper}) is
fulfilled and the eigenvalues $\mu_\k(\la,m)$ (\ref{mu-lam}) are all
distinct.
\end{definition}

So, in an $m$-admissible representation $\opi_\la$ we have
\be
\opi_\la(\Tr L_{(m)}^s)=\sum_{|\k|=m} \mu_\k(\la,m)^s d_\k(\la,m)
\,{\rm Id}_{V_\lambda^{\!*}}
\label{trlmk}
\ee
with some multiplicities $d_\k(\la,m)$. The identity operator in the
right hand side appears since $\Tr L_{(m)}^s\in Z(U(\gg_\hbar))$ and
$\opi_\la$ is an irreducible representation.

Now, on the one hand, the multiplicities $d_\k(\la,m)$ should be
specializations of  $d_\k(m)$ (\ref{razmer}) at
$\mu_i = \mu_i(\lambda)$. On the other hand,  $d_\k(\la,m)$ can be
computed independently by means of  the split Casimir element. To
prove the proposition, we have to show that these two ways of
computation give the same results.

So, we proceed in the way analogous to proof of proposition
\ref{prop:3}. We apply the representation $\opi_\lambda$ to the
matrix $L^s_{(m)}$ passing thereby to $L^s_{(\lambda,m)}$. Then, the
transposed matrix $(L^s_{(\lambda,m)})^t$ is nothing but the matrix
of $\Caslm^s$, i.e., the $s$-th power of the linear operator
(\ref{cas-lr}) constructed via the split Casimir element (for
detail, see \cite{GLS2}).

Applying the transposition operator to the CH identities $\chc_{(m)}
(L_{(\lambda,m)})=0$ where $\chi$ is fixed by the values of $\mu_i$
(\ref{mu}), we conclude that the set of eigenvalues of the matrix
$L_{(\lambda,m)}$ coincides with that of the operator $\Caslm$.

Since the operator $\Caslm$ is scalar on irreducible components
$V_{\lambda^{\!*}+\k}$, we get
\be
\Trr\,\Caslm^s=\sum_{|\k|=m} \mu_\k(\la,m)^s\dim V_{\la^{\!*}+\k}.
\label{altern}
\ee
Here $\Trr=\Tr\ot\Tr:\End(\vl)\ot\End(\vm)\to\K$. Taking the traces
in (\ref{trlmk}) and comparing the result with (\ref{altern}) we
conclude that
\be
d_\k(\la,m) = \frac{\dim V_{\la^{\!*}+\k}}{\dim V_{\la^{\!*}}}.
\label{dk}
\ee

At last, using the well known Frobenius formula
$$
\dim V_{\la} = \prod_{i<j} {(\la_i-\la_j-i+j)\over (j-i)}
$$
and taking into account (\ref{mu}), we find that $d_\k(\la,m)$ are
indeed given by (\ref{razmer}) with the substitution
$\mu_i = \mu_i(\la)$ and $\h=1$. The passage to a generic $\h\not=0$
is standard.

So, formulae (\ref{razmer}) and (\ref{princip}) have been proved for
any $m$-admissible representation. But since the family of
$m$-admissible signatures $\la$ is a large enough (actually
infinite) set, then by virtue of reasons discussed at the beginning
of the proof we conclude that  (\ref{princip}) is true in general.
\hfill\rule{6.5pt}{6.5pt}

\medskip
If $m=1$, then (\ref{princip}) reduces to
\be
\Tr L^s=\sum_{j=1}^n \mu_j^s \prod_{i\not=j}
\frac{\mu_j-\mu_i-\h}{\mu_j-\mu_i}.
\label{reduc}
\ee
As we have said above, at $\h=1$ and $\mu_i=\mu_i(\la)$ this formula
is equivalent to that from \cite{BR}.

Formula (\ref{reduc}) together with (\ref{elsym}) can be treated as
an $n$-parametric resolution of the Newton relations in the algebra
$U(\gh)$. Indeed, the right hand side of (\ref{reduc}) is a
symmetric polynomial in $\mu_i$, and therefore it can be expressed
as a polynomial in elementary symmetric functions in $\mu_i$. Upon
replacing these functions by $\sigma_k$, we find explicit relations
between two central families $\{\sigma_k\}$ and $\{s_k=\Tr L^k\}$ in
the algebra $U(\gh)$. The variables $\mu_i$ can be thought of as the
elements of the algebraic closure of the center of $U(\gh)$.

In the same sense we consider (\ref{princip}) as the higher order
counterparts of the Newton relations. In section \ref{sec:6} we
present the $q$-analogs of these relations.

\section{Reflection Equation Algebra: definition and basic
properties}
\label{sec:rea-def}

This section is devoted to the detailed description of the Hecke
symmetries $R$ and corresponding reflection equation algebras which
were shortly outlined in the Introduction.

Let $V$ be a vector space over $\K$,  $\dim\, V=n$ and let
$R \in \End(\vv)$ be an endomorphism. We call $R$ {\em braiding} if
it satisfies the {\it Yang-Baxter equation}
\be
R_{12}R_{23}R_{12} = R_{23}R_{12}R_{23},
\label{YBE}
\ee
where $R_{12}=R\ot {\rm id}_V$ and $R_{23}={\rm id}_V\ot R$ are
treated as elements of $\End(V^{\ot 3})$. On fixing a basis
$\{x_i\}\in V$, $1\le i \le n$, in the space $V$, we can realize the
endomorphism $R$ as a numerical $n^2\times n^2$ matrix $R$ for which
we use the same notation and call it {\em $R$-matrix}:
\be
(x_i\ot x_j)\triangleleft R=R^{\,kl}_{ij}\, x_k\ot x_l.
\label{R-matr}
\ee

In what follows we shall deal with a special case of {\it Hecke
type} $R$-matrices satisfying the four additional conditions listed
below under the items {\bf C1) -- C4)}.

\medskip

\noindent{\bf C1)}\
First of all, the $R$-matrix should obey the {\it Hecke condition}
\be
(R-qI)(R+q^{-1}I) = 0,
\label{Hecke}
\ee
where $q$ is a fixed nonzero number from the ground field $\Bbb K$
with the only constraint
\be
q^m\not=1, \quad\forall\,m\in{\Bbb N}.
\label{generic}
\ee
As a consequence, the $q$-analogs of all integers are nonzero
\be
m_q\equiv\frac{q^m-q^{-m}}{q-q^{-1}}\not=0,\quad \forall\, m\in
{\Bbb Z}.
\label{q-num}
\ee
In some cases it proves to be convenient to consider $q$ as a formal
parameter and extend $\Bbb K$ to the field of rational functions in
the indeterminate $q$. We shall always bear in mind this extension
when considering the classical limit $q\rightarrow 1$.

\medskip

\noindent{\bf C2)}\
To formulate the further restriction on $R$, we need to recall some
properties of the Hecke algebras of $A_{m-1}$ series and to describe
their relation to the Hecke type $R$-matrices.

Fix a nonzero number $q\in \Bbb K$. The {\it Hecke algebra of
$A_{m-1}$ series} ($m\ge 2$) is an associative algebra $H_m(q)$ over
the field $\Bbb K$ generated by the unit element $1_{H}$ and $m-1$
generators $\sigma_i$ subject to the following relations:
$$
\left.
\begin{array}{l}
\sigma_i\sigma_{i+1}\sigma_i = \sigma_{i+1}\sigma_i\sigma_{i+1}\\
\rule{0pt}{5mm}
\sigma_i\sigma_j = \sigma_j\sigma_i\qquad\qquad {\rm if\ }|i-j|\ge
2\\
\rule{0pt}{5mm}
(\sigma_i-q\,1_{H})(\sigma_i+q^{-1}\,1_{H}) = 0\\
\end{array}
\right\}\quad  i=1,2,\dots m-1.
$$

If the parameter $q$ satisfies (\ref{generic}), then for any
positive integer $m$ the Hecke algebra $H_m(q)$ is isomorphic to the
group algebra ${\Bbb K}S_m$ of the $m$-th order permutation group
$S_m$. As a consequence, $H_m(q)$ is isomorphic to the following
direct sum of matrix algebras
\be
H_m(q)\cong \bigoplus_{\lambda\vdash m}{\rm Mat}_{d_\lambda}(\Bbb
K)
\label{H-iso-Mat}
\ee
where the summation goes over all possible ordered partitions
$\lambda$ of the integer $m$. The parameter $d_\lambda$ is equal to
the number of all standard Young tableaux $\lambda(a)$ which can be
constructed for the given partition $\lambda$.

As is known, the associative $\Bbb K$-algebra ${\rm Mat}_k(\Bbb K)$
possesses the linear basis of $k^2$ generators (matrix units)
$E_{ab}$ $1\le a,b \le k$ with the multiplication law
$$
E_{ab}E_{cd} = \delta_{bc}E_{ad}.
$$
Due to isomorphism (\ref{H-iso-Mat}) in the Hecke algebra $H_m(q)$
one can choose a system of generators
${\cal Y}^\lambda_{ab}(\sigma) \in H_m(q)$, $\lambda\vdash m$,
$1\le a,b \le d_\lambda$, which form a linear basis in $H_m(q)$ and
satisfy the following multiplication rule
\be
{\cal Y}^\lambda_{ab}(\sigma){\cal Y}^\mu_{cd}(\sigma) =
\delta^{\lambda\mu}\delta_{bc} {\cal Y}^\lambda_{ad}(\sigma).
\label{Y-Y}
\ee
The diagonal "matrix units" ${\cal Y}^\lambda_{aa}(\sigma)$ are the
primitive idempotents of the Hecke algebra. Evidently, they are in
one-to-one correspondence with the set of standard Young tableaux
$\lambda(a)$ constructed for the partition $\lambda$. Below we shall
use the compact notation
${\cal Y}^\lambda_{aa}\equiv{\cal Y}_{\lambda(a)}$ for the primitive
idempotents. The idempotent, corresponding to the partition $(1^m)$
(the one-column Young diagram) will be called the
$q$-antisymmetrizer and denoted ${\cal A}^{(m)}(\sigma)$.

Given a Hecke type $R$-matrix, we can construct a {\em local
representation} of  $H_m(q)$ in $V^{\otimes m}$ by the following
rule
\be
\sigma_i\mapsto \rho_R(\sigma_i)= R_{ii+1} = I^{\otimes(i-1)}
\otimes
R\otimes I^{\otimes(m-i-1)}\in {\rm End}(V^{\otimes m}).
\label{loc-rep}
\ee
In the local representation (\ref{loc-rep}) the idempotents
${\cal Y}_{\lambda(a)}(\sigma)$ are realized as some projection
operators in $V^{\otimes m}$. With respect to the action of these
projectors the space $V^{\otimes m}$ splits into the direct sum of
subspaces $V_{\lambda(a)}$:
\be
V^{\otimes m} = \bigoplus_{\lambda\vdash m} \bigoplus_{a=1}
^{d_\lambda}
V_{\lambda(a)},\quad V_{\lambda(a)} = Y_{\lambda(a)}(R)
\triangleright
V^{\otimes m}. \label{razl-prostr}
\ee
The projector $Y_{\lambda(a)}(R) = \rho_R({\cal Y}_{\lambda(a)})$ is
given by a polynomial in matrices $R_{ii+1}$. For a detailed
treatment of these questions, explicit formulae for $q$-projectors
and the extensive list of original papers, see \cite{OgP}.

So, we shall assume the Hecke symmetry in question to be even. By
definition this means that  there exists an integer $p>0$ such that
the image of the $q$-an\-ti\-sym\-met\-ri\-zer
${\cal A}^{(p+1)}(\sigma)$ in the local $R$-matrix representation
$\rho_R$ identically vanishes while the image of the
$q$-antisymmetrizer ${\cal A}^{(p)}(\sigma)$ is a projector of the
unit rank in the space $V^{\otimes m}$ (for any $m>p$)
\be
\exists\,p\in{\Bbb N}:\quad
\left\{
\begin{array}{l}
{\cal A}^{(p+1)}(\sigma)\mapsto A^{(p+1)}(R)\equiv 0,\\
\rule{0pt}{6mm} {\cal A}^{(p)}(\sigma)\mapsto A^{(p)}(R),
\quad {\rm rank}\,A^{(p)}(R) = 1.
\end{array}
\right.
\label{f-rank}
\ee
Such a number $p$ will be called {\it the symmetry rank} of the
matrix $R$. For example, the symmetry rank of the $R$-matrix
connected with the quantum universal enveloping algebra $U_q(sl(n))$
is equal to $n$. Examples of $n^2\times n^2$ $R$-matrices with $p<n$
(for $n\geq 3$) can be found in \cite{G}.

\medskip

\noindent{\bf C3)}
Besides, we assume the $R$-matrix to be {\it skew-invertible}; that
is there exists an $n^2\times n^2$ matrix $\Psi$ with the property
$$
\sum_{a,b}R_{ia}^{\;jb}\Psi_{bk}^{\;as} = \delta_i^{\,s}
\delta_k^{\,j} =
\sum_{a,b} \Psi_{ia}^{\;jb}R_{bk}^{\;as}.
$$
In the compact notations the above formula reads
\be
\Tr_{(2)}R_{12}\Psi_{23} = P_{13} = \Tr_{(2)}\Psi_{12}R_{23},
\label{closed}
\ee
where the symbol $\Tr_{(2)}$ means applying the trace in the second
space and $P$ is the permutation matrix. Note, that this condition
does not depend on a choice of  basis, therefore we can treat $\Psi$
as an endomorphism as well.

Let now $B$ and $C$ be two endomorphisms of the space V represented
by the following $n\times n$ matrices
\be
B = \Tr_{(2)}\Psi_{21}, \quad C = \Tr_{(2)}\Psi_{12},
\label{BC}
\ee
where $\Psi$ is defined in (\ref{closed}). If the $R$-matrix has the
symmetry rank $p$, we have
\be
B\cdot C = \frac{1}{q^{2p}}\,I
\label{nonsing}
\ee
(cf. \cite{G}). This implies that the matrices $B$ and $C$ are
invertible. Moreover, we have
\be
\Tr B=\Tr C=\frac{p_q}{q^p}.
\label{BC-norm}
\ee
Note that {\bf C3)} is in fact a consequence of {\bf C2)} since any
even Hecke symmetry is skew-invertible automatically.
\medskip

\noindent{\bf C4)}\
As the last requirement on the $R$-matrix we shall assume, that
the space $V^*$ dual to $V$ can be identified with the $(p-1)$-th
wedge $q$-power of $V$
$$
V^* = \wedge^{p-1}_q V.
$$
Explicitly this identification is constructed as follows. Since
$A^{(p)}(R)$ is the unit rank projector then in the basis
$x_{i_1}\otimes\dots \otimes x_{i_p}$ (see (\ref{R-matr})) its
matrix can be represented in the form
$$
{A^{(p)}}_{i_1i_2\dots i_p}^{\;\;j_1j_2\dots j_p}
 = u_{i_1i_2\dots i_p}v^{j_1j_2\dots j_p}
$$
with some structure tensors $u$ and $v$. One can show that the
vectors
\be
x^i = v^{ia_2\dots a_p}x_{a_2}\otimes \dots \otimes x_{i_p}
\label{dual-basis}
\ee
are linear independent and by definition form the basis of the space
$\wedge^{p-1}_q V\subset V^{\otimes (p-1)}$. We shall not explicitly
describe the property of $R$ which allows us to identify $V^*$ and
$\wedge^{p-1}_q V$. The thorough treatment of this problem is
presented in \cite{GLS1}.

In the particular case $p=2$ which will be studied in detail below
the components of structure tensor $v^{ij}$ form a nondegenerated
matrix \cite{G}. Therefore, as follows from (\ref{dual-basis}), the
space $V^*$ is isomorphic to $V$ itself
\be
V^*\cong V\quad {\rm at}\quad p=2.
\label{v-dua-v}
\ee

\medskip

As the next step, we define the reflection equation algebra.
Consider a unital associative $\K$-algebra $\lhq$ generated by $n^2$
elements $l_i^{\,j}$,  $1\le i,j\le n$ satisfying the following
relations
\be
R L_1R L_1 - L_1R L_1R = \hbar(RL_1
- L_1R), \quad L_1\equiv  L\otimes I,
\label{RE}
\ee
where $\hbar$ is a numerical parameter and $ L =\|l_i^{\,j}\|$ is a
matrix composed of $l_i^{\,j}$. The algebra $\lhq$ is said to be the
{\it modified reflection equation algebra} (mREA). In the particular
case $\h= 0$ this algebra will be called (non-modified) REA and
denoted $\lq$.

\begin{remark}\label{rem-1}
{\rm Similarly to the algebra $U(\gh)$ the mREA corresponding to
different {\it nonzero} parameters $\hbar$ are isomorphic --- one
can easily pass from one $\hbar\not=0$ to another $\hbar'\not =0$ by
a trivial renormalization of generators.

Moreover, at $q\not=\pm 1$ the generators $l_i^{\,j}$ of the
algebra \lqh{} is connected with the generators $\hat l_i^{\,j}$
of $\lq$ via a linear shift by the unit element ${\rm id}_{\cal L}$
\be
l_i^{\,j} = \hat l_i^{\,j} + \frac{\hbar}{\zeta} \,\delta_i^{\,j}\,
{\rm id}_{\cal L},\qquad \zeta = q-q^{-1},
\label{shift}
\ee
and therefore these two kinds of REA are isomorphic, too.
Nevertheless, their classical limits are different, since the above
isomorphism is broken as $q\rightarrow 1$.}
\end{remark}

Let us now define an important  map $\Tr_R:{\rm Mat}_n(\lqh)
\rightarrow \lqh$
\be
\Tr_R(X)\stackrel{\mbox{\tiny def}}{=} \Tr (C\cdot X),\quad X\in
{\rm
Mat}_n(\lqh).
\label{q-sled}
\ee
For the $U_q(sl(n))$ $R$-matrix, such a map is called the
{\it quantum trace} \cite{FRT} and is often denoted $\Tr_q$. Note
that in the sequel we apply the trace (\ref{q-sled}) only to the
matrices from the space ${\rm Mat}_n(\lqh)$ which are in some sense
invariant.

Similarly to the classical case, a generating set for the center of
the algebra \lqh{} consists of the elements:
\be
s_0\equiv {\rm id}_{\cal L},\quad
s_k = q\,\Tr_RL^k,\quad 1\le k\le p,
\ee
where the factor $q$ is chosen for the future convenience. These
elements become scalar operators on each irreducible module over
the mREA (see section \ref{sec:rea-br}).

One of the basic properties of the REA (modified or not) consists in
the following. In this algebra, similarly to $U(su(n))$ (and to some
other Lie algebras, cf. \cite{Go}), one can find a series of
the Cayley-Hamilton identities. Let us discuss the first of them,
which will be called the {\it basic CH identity}.

As was shown in \cite{GPS}, the matrix $\hat L=\|\hat l_i^{\,j}\|$
of generators of REA satisfies the CH identity of the form
\be
\sum_{k=0}^p(-\hat
L)^{p-k}\sigma_k(\hat L) = 0,\quad \sigma_0(\hat L) \equiv
{\rm id}_{\cal L},\quad \hat L^0\equiv I
\label{CH}
\ee
where $\{\sigma_k\}$ is another set of generators of the center of
REA connected with the set $\{s_k\}$ by means of the $q$-Newton
relations (\ref{NI}) and $p$ is the symmetry rank of the
corresponding $R$.

Applying the shift (\ref{shift}) it is possible to get analogous
relations for the matrix $L$ with entries belonging to the
corresponding mREA:
$$
\sum_{k=0}^p(-L)^{p-k}\sigma_k(L) = 0,\qquad \sigma_0(L)
\equiv {\rm id}_{\cal L}.
$$
As above the coefficients $\sigma_k(L)$ are central element of the
algebra $\lhq$.

An explicit form of these coefficients was obtained in \cite{GS1}.
However, this form is somewhat cumbersome. Here we only present an
example of such a CH identity for the mREA algebra related to a
Hecke symmetry with the symmetry rank $p=2$:
\be
L^2 - \Big(q\Tr_R L +\frac{\hbar}{q}\Big)\,L+\Big(\frac{q^2}{2_q}
\Big(q(\Tr_RL)^2 - \Tr_RL^2\Big)+ \hbar\,\frac{q}{2_q}\,
\Tr_R L \Big)\,I = 0.
\ee
In this formula the coefficients $\sigma_i(L)$ are expressed via the
quantities $\Tr_R(L^k)$. In what follows we will find some inverse
relations expressing the latter quantities via $\sigma_i(L)$ but in
a parametric form.

By analogy with the case of fuzzy orbit (see definition
\ref{def:fuz-orb} and formula (\ref{gen})) we define a character
$\chi: Z(\lqh) \to {\Bbb K}$ of the center of the mREA \lqh{} by
fixing its values on the central elements $\sigma_k$
\be
\chi(\sigma_k(L)) = \alpha_k = \sum_{1\le i_1<\dots <
i_k\le p}\mu_{i_1} \dots \mu_{i_k},
\label{q-gen}
\ee
where the numbers $\mu_i$ $1\le i\le p$ are assumed to be all
distinct.

Then we define a {\it 1-generic NC orbit} $\lhqc$ as the quotient
of the mREA (\ref{RE}) modulo the two sided ideal ${\cal I}^\chi$
\be
\lhqc = \lqh/{\cal I}^\chi
\label{q-orb}
\ee
where ${\cal I}^\chi$ is generated by the set of relations
\be
\sigma_k(L)-\alpha_k, \quad {\rm for } \quad 1\le k
\le p.
\label{ideal}
\ee
Finally, we will compute the quantities $\Tr_R L^k$ restricted on
the NC orbit $\lhqc$ in terms of $\mu_i$.

Following the pattern of the case considered in section
\ref{sec:fuzzy} we will also compute the quantities
$\Tr_R L^k_{(m)}$ for some higher extensions $L_{(m)}$ of the
matrix $L$. The crucial role in this computing is played by a split
Casimir element which is defined for the mREA in the following way
\be
{\bf Cas} = q^{2p}\,l_i^{\,k}\otimes l_k^{\,j}C_j^{\,i}\;\in
\;\lqh\otimes \lqh.
\label{quant-Cas}
\ee
As in the classical case considered in section \ref{sec:fuzzy}, this
element allows us to introduce the aforementioned extensions
$L_{(m)}$ of the matrix $L$  and to construct the  higher CH
identities for them. Also note, that in the $\Uq$ case at the limit
$q=1$ we get just the split Casimir element  (\ref{cas-kl})
described above.

\section{mREA as a braided enveloping algebra}
\label{sec:rea-br}

In this section we give a short review of the representation theory
of REA which is necessary for the subsequent sections. Besides, we
adduce some arguments which allow us to consider the mREA as a
"braided" analog of the enveloping algebra.

First of these arguments originates from the consideration of the
mREA for an {\it involutive} $R$-matrix: $R^2=I$. This is a
particular case of the Hecke condition (\ref{Hecke}) corresponding
to $q=1$. The mREA (\ref{RE}) with involutive $R$-matrix turns out
to be the enveloping algebra of some generalized Lie algebra.

The notions of generalized Lie algebra and its enveloping algebra
were introduced in \cite{G} (see also the references therein). Given
an involutive $R$-matrix (treated as an endomorphism of
$V^{\otimes 2}$), the generalized enveloping algebra is defined as
the quotient of the free tensor algebra of the space $\End(V)$ over
the two-sided ideal generated by the relations
$$
X\ot Y-R_{\End(V)}(X\ot Y)=\circ X\otimes Y-\circ
R_{\End(V)}(X\ot Y),\quad  X,Y\in \End(V),
$$
where $\circ:\End(V)^{\ot 2}\to \End(V)$ is the usual product in the
space of endomorphisms and $R_{\End(V)}$ is the extension of the
initial braiding $R$ to the space $\End(V)$ (which is well defined
as $R$ is skew-invertible).

If we realize $\End(V)$ as the space of left endomorphisms, fix a
natural basis $\{h_i^{\,j}\}$ such that $h_i^{\,j}\circ h_k^{\,l} =
\delta_k^j h_i^{\,l}$ and compute $R_{\End(V)}$ in this basis, we
recover the enveloping algebra of the generalized Lie algebra in
terms of generators $h_i^{\,j}$. Note, that $h_i^{\,j}$ can be
identified with $x_i\ot x^j$ where $\{x^j\}$ is the basis of the
left dual space to $V$, see \cite{G} for detail.

The point is that in $\End(V)$ we can choose another basis
$\{l_i^{\,j}\}$ such that
$l_i^{\,j}\circ l_k^{\,l} = l_i^{\,l}\, B_k^j$, where $B_i^j$ is
the matrix element of the endomorphism $B$ (\ref{BC}) written in the
basis $\{x_i\}$. Being expressed in terms of the new generators
$l_i^{\,j}$, the enveloping algebra in question turns into the mREA
with $\h=1$. The elements $l_i^{\,j}$ can be identified with
$x_i\ot{^j}x$ where $\{{^j}x\}$ is the  basis of the right dual
space to $V$.) The details are left to the reader.

Let us also mention that in the $\Uq$ case (in contrast with quantum
groups corresponding to other simple Lie algebras) the mREA is a
two parameter deformation of the symmetric algebra of $gl(n)$ (a
close treatment of this algebra is given in \cite{IP}). The
corresponding Poisson structure is the aforementioned pencil which,
in fact, is well defined on the whole $gl(n)^*$.

The most important property of the mREA which enable us to treat
this algebra as a braided analog of the enveloping algebra is that
the category of its equivariant finite dimensional representations
is close to that of $U(gl(p))$--Mod where $p$ is the symmetry rank
of the $R$-matrix.

Let us first consider the quotient of the mREA \lqh{} over the ideal
generated by the relation $\Tr_RL = 0$. We denote this quotient
algebra as ${\cal SL}_{q,\h}$. This is an analog of the $U(sl(p))$
subalgebra in the $U(gl(p))$. The category of finite dimensional
completely reducible modules over the ${\cal SL}_{q,\h}$ is the
{\it quasitensor Schur-Weyl category} introduced in \cite{GLS1}.
Its simple objects (irreducible modules) are labelled by the
partitions $\lambda$ whose height (the number of nonzero parts) are
not greater than $p-1$. Besides, the Grothendieck ring of the
Schur-Weyl category is isomorphic to that of the category
$U(sl(p))$--Mod.

The representations of \lqh{} are labelled by the partitions
$\lambda$ whose height is not greater than $p$ and a number
$z\in {\Bbb K}$ which is analog of the shift (\ref{ed-sdv}) of the
$U(gl(p))$ representations (see remark \ref{rem:14} below). In full
analogy with the classical case discussed in remark \ref{rem:3}, the
vector spaces $V_{\lambda,z}$ and $V_{\hat \lambda}$ are isomorphic
as ${\cal SL}_{q,\h}$ modules, $\hat\lambda$ being constructed from
$\lambda$ in accordance with (\ref{lambda-hat}).

The map ${\rm Tr}_R$ defined in (\ref{q-sled}) is closely related to
the {\em categorical trace}
$\mbox{\bf\sf Tr}_{V_\lambda}:\End(V_\lambda)\to{\Bbb K}$. The
categorical trace is a morphism of the Schur-Weyl category which
plays the same role as the usual trace does in the category
$U(sl(p))$--Mod. In particular, the categorical trace allows one
to define the notion of the {\it $q$-dimension}
\be
\dim_qV_\lambda = \mbox{\bf\sf Tr}_{V_\lambda}
({\rm id}_{V_\lambda}),
\label{def:q-dim}
\ee
which is a multiplicative-additive functional on the Grothendieck
ring  of the Schur-Weyl category. Namely we have
\be
\begin{array}{l}
\dim_q(V_\lambda\otimes V_\mu)=\dim_q(V_\lambda)\,\dim_q(V_\mu),\\
\rule{0pt}{5mm}\dim_q(V_\lambda\oplus V_\mu) = \dim_q(V_\lambda)+
\dim_q(V_\mu)
\end{array}
\label{ma-func}
\ee
and the $q$-dimensions of the isomorphic spaces are equal to each
other (see \cite{GLS1} for detail).

Now let us give a short review of some facts from the representation
theory of the mREA. Their detailed description can be found in
\cite{S}.

The mREA \lqh{} possesses a profound representation theory. We shall
confine ourselves to considering the {\it finite dimensional,
completely reducible, equivariant} modules over $\lqh$  (the term
"equivariant" will be explained at the end of the section). Besides,
we can set $\hbar = 1$ by virtue of the isomorphisms mentioned in
Remark~\ref{rem-1}.

First of all, we define the so-called {\it left fundamental module
of B type}. Let $V$ be an $n$-dimensional vector space with a fixed
basis $\{x_i\}$ $1\le i\le n$. Putting $\h=1$, we consider the
homomorphism $\pi:\loq\rightarrow {\rm End}_l(V)$ defined as follows
\be
\pi(l_i^{\,j})\triangleright x_k = x_iB_k^{\,j},
\label{B-rep}
\ee
where the matrix $B$ is introduced in (\ref{BC}). Since $B$ is
invertible (see (\ref{nonsing})) the representation $\pi$ is
irreducible.

The tensor power $V^{\otimes m}$, $m\in {\Bbb N}$, can be endowed
with the structure of a (reducible) mREA module. The corresponding
homomorphism $\rho_m:\loq\rightarrow {\rm End}_l(V^{\otimes m})$
is of the form
\be
\rho_m(l_i^{\;j}) = \pi_1(l_i^{\;j}) + {\cal R}_{12}^{-1}
\pi_1(l_i^{\;j}) {\cal R}_{12}^{-1}  + \dots +
{\cal R}_{m-1,m}^{-1}\dots {\cal R}_{12}^{-1}\pi_1(l_i^{\;j})
{\cal R}_{12}^{-1}\dots {\cal R}_{m-1,m}^{-1},
\label{m-rep}
\ee
where
$$
\pi_1 = \pi \otimes I^{\otimes (m-1)}, \quad
{\cal R}_{k,k+1} = I^{k-1}\otimes {\cal R} \otimes
I^{m-k-1}\quad 1\le k\le m-1.
$$
Here ${\cal R}$ is an automorphism of $V^{\otimes 2}$ connected with
the matrix $R$ of the Hecke symmetry by the following definition
\be
{\cal R}\triangleright(x_i\otimes x_j) = \sum_{r,s}R_{ij}^{\;kl}
x_k\otimes x_l.
\label{R}
\ee

The representation (\ref{m-rep}) is reducible. In accordance with
(\ref{razl-prostr}) it decomposes into the direct sum of mREA
submodules $V_{\lambda(a)}$, where $\lambda$ is an ordered partition
of $m$. The representations $\pi_{\lambda(a)}$ are extracted from
$\rho_m$ by the action of the corresponding projectors
$Y_{\lambda(a)}(R)$.

We write down the explicit form of the representation $\pi_{(m)}$,
corresponding to the partition $(m)$, since it will play an
important role in what follows. The subspace
$V_{(m)}\subset V^{\otimes m}$ is an image of the
{\it $q$-symmetrizer} $S^{(m)}({\cal R})$ whose matrix is
iteratively defined as follows (see \cite{G})
\be
S^{(1)}\equiv I,\qquad
S^{(m)}_{12\dots m} = \frac{1}{m_q}\,S^{(m-1)}_{2\dots m}(q^{1-m}
I+(m-1)_q
R_{12})S^{(m-1)}_{2\dots m}.
\label{q-symm-S}
\ee
The following proposition holds true \cite{S}.

\begin{proposition} \label{prop:pim}
Consider an arbitrary tensor power $V^{\otimes m}$ of the left
fundamental module $V$. Its $q$-symmetric subspace $V_{(m)}$ is a
left $\loq$-submodule. On the generators of the mREA the
homomorphism $\pi_{(m)}:\loq \rightarrow {\rm End}_l(V_{(m)})$ reads
as follows:
\be
\pi_{(m)}(l_i^{\,j}) = q^{1-m}m_q\,S^{(m)}({\cal R})
\Big[\pi(l_i^{\,j})
\otimes I^{\otimes (m-1)}\Big]S^{(m)}({\cal R}),
\ee
where $\pi$ and ${\cal R}$ are defined in (\ref{B-rep}) and
(\ref{R}) respectively.
\end{proposition}

Given the representation (\ref{B-rep}), we can realize the matrix
$L$ satisfying the commutation relations (\ref{RE}) as an image of
the split Casimir element {\bf Cas} (\ref{quant-Cas}) under the map
$$
{\rm id}\otimes \pi:\; \lqh\otimes\loq\rightarrow \lqh\otimes
{\rm Mat}_n(\Bbb K).
$$
(Note that in the second factor we put $\h=1$.) Indeed, taking
(\ref{nonsing}) into account one easily gets
\be
L^t = ({\rm id}\otimes \pi)({\bf Cas}),\quad L\in
{\rm Mat}_n(\lqh).
\label{L-real}
\ee
Hereafter $L^t$ stands for the matrix transposed to $L$.

By analogy with the  $U(\gg)$ case considered in section
\ref{sec:fuzzy} we also introduce the symmetric matrix $L_{(m)}$ as
the following image of ${\bf Cas}$
\be
L_{(m)}^t = ({\rm id}\otimes \pi_{(m)})({\bf Cas}),
\quad L_{(m)}\in {\rm Mat}_{n_m}(\lqh)
\label{L-symm}
\ee
where $n_m = {\dim}\,V_{(m)}$.

In addition to the left $\loq$-modules we also need the right
ones. For a generic $p\geq 2$ such a representation can be defined
in the dual space $V^*$. By using the method of the paper \cite{S}
we can extend this representation to the tensor power
$(V^*)^{\ot m}$ and decompose it into the direct sum of submodules
associated with the projectors $Y_\lambda$. However, if $p=2$ the
space $V^*$ can be identified with $V$ itself via a categorical
pairing arising from the projector of the space $V^{\ot 2}$ onto its
skew-symmetric component (see \ref{v-dua-v}). The spaces
$V_{(k)}^*$ and $V_{(k)}$ can also be identified via such a pairing
(see \cite{GLS1} for detail). This pairing allows us to equip the
space $V_{(k)}$ with a structure of the right $\loq$-module. More
precisely, the following proposition takes place.

\begin{proposition}\label{prop:r-m}
Let the symmetry rank of the Hecke symmetry $R$ is $p=2$. Consider
an arbitrary tensor power $V^{\otimes m}$ of the left fundamental
module $V$. Its $q$-symmetric subspace $V_{(m)}$ is a right
$\loq$-submodule. On the generators of the mREA the homomorphism
$\overline\pi_{(m)}:\loq \rightarrow {\rm End}_r(V_{(m)})$ reads as
follows:
\be
\overline\pi_{(m)}(l_i^{\,j}) = q^{1-m}m_q\,S^{(m)}(R)
\Big[I^{\otimes (m-1)}
\otimes \overline\pi(l_i^{\,j})\Big]S^{(m)}(R),
\ee
where the homomorphism
$\overline\pi:\loq \rightarrow {\rm End}_r(V)$ has the form
\be
x_k\triangleleft \overline\pi(l_i^{\,j})=
\frac{2_q}{q^2} {A^{(2)}}_{ki}^{\;sj}x_s,\qquad
A^{(2)}\equiv \frac{1}{2_q}\,(qI - R),
\label{rep:pi-r}
\ee
and the right action of $R$ is defined in (\ref{R-matr}).
\end{proposition}

In what follows  we shall also need the representations of the
generators $\hat l_i^{\,j}$ of the REA connected with those of mREA
by shift (\ref{shift}). A simple calculation proves the following
corollary of the proposition \ref{prop:r-m}.

\begin{corollary}
\label{prop:rea-r}
The right representation of the generators $\hat l_i^{\,j}$ of the
REA in the $q$-symmetric component $V_{(m)}\subset V^{\otimes m}$ is
given by the homomorphism
\be
\overline \pi_{(m)}(\hat l_i^{\,j}) =
q^{1-m}m_q\,S^{(m)}(R)\Big[I^{\otimes (m-1)} \otimes
\overline\pi(\hat l_i^{\,j})\Big]S^{(m)}(R),
\label{rep:pim-r}
\ee
where
$$
x_k\triangleleft \overline\pi(\hat l_i^{\,j}) =
\Phi_{ki}^{\;sj}x_s ,\quad \Phi \equiv  q^{1-m}m_q\,I -
\zeta\, \frac{2_q}{q^2}\,A^{(2)}.
$$
\end{corollary}

\begin{remark}\label{rem:14}
{\rm
Up to now, the representations of mREA were labelled by partitions
$\lambda$, whereas in general the finite dimensional representations
of $U(gl(n))$ are labelled by {\it signatures} $\lambda$
(\ref{signature}). What is an analog of these representations in the
mREA case? To answer the question, observe the following.

As was mentioned in remark \ref{rem:3}, with any finite dimensional
$U(gl(n))$ representation $\pi_\lambda$ labelled by a signature
$\lambda$ we can associate the representation $\pi_{\hat \lambda}$
where the partition $\hat \lambda$ is connected with the signature
by relation (\ref{lambda-hat}). These representations are connected
by the unit ope\-ra\-tor shift (\ref{ed-sdv}) and the
corresponding modules $V_\lambda$ and $V_{\hat \lambda}$ are
isomorphic as $U(sl(n))$-modules and therefore as vector spaces.

For the mREA representations there exists a transformation analogous
to (\ref{ed-sdv}). Namely, a simple calculation shows that if
$\pi_\lambda$ is a representation of mREA in the space $V_\lambda$
where $\lambda$ is a partition then the operators
\be
\pi_\lambda^z(l_i^{\,j}) = z\,\pi_\lambda(l_i^{\,j}) +
\delta_i^{\,j}\, \frac{1-z}{\zeta}\,{\rm id}_{V_\lambda},
\quad z\in{\Bbb K}\backslash 0
\label{pi-z}
\ee
also realize an mREA representation in the same space $V_\lambda$.
It can be shown that the representation of the subalgebra
${\cal SL}_{\hbar,q}$ which is a quotient of $\lqh$ over the ideal
generated by the relation ${\rm T}_RL = 0$ does not change under the
shift (\ref{pi-z}). This subalgebra is an analog of $U(sl(n))$ in
the classical case. However, in contrast with the classical case
discussed in remark \ref{rem:3}, we cannot put any signature into
correspondence to $\pi_\lambda^z$ since our approach to the
representation theory of the mREA is not based on the technique of
the highest weight vectors.}
\end{remark}

Now we explain the meaning of the statement that the proposed
presentation theory of mREA is equivariant.

To any $R$-matrix we can assign an associative bialgebra $\cal T$
generated by the elements $t_i^{\,j}$ subject to the following
commutation relations
$$
R\, T_1 T_2 = T_1 T_2R, \quad {\rm where}\quad T=
\|t_i^{\,j}\|,\quad T_1 = T\otimes I,\quad T_2 = I\otimes T.
$$
If the $R$-matrix is skew-invertible (see (\ref{closed})), then
the bialgebra structure can be extended to the Hopf algebra
one\footnote{In order to get such a structure it suffices to
formally invert the quantum determinant in the extended algebra
$\cal T$. It can be done by an appropriate localization, cf.
\cite{G}.}. When $R$ is the image of the universal $U_q(sl(n))$
$R$-matrix (in the fundamental vector representation), the Hopf
algebra $\cal T$ is the well known quantization of the algebra of
regular functions on the group $GL(n)$ \cite{FRT}. Since the REA can
always be endowed with the structure of the left adjoint comodule
over $\cal T$
$$
l_i^{\,j}\mapsto \sum_{r,s}t_i^{\,r}S(t_s^{\,j})\otimes l_r^{\,s},
$$
$S(t_i^{\,j})$ being the antipode of $t_i^{\,j}$, then REA is also a
module over the dual Hopf algebra $\cal T^*$ (in the $U_q(sl(n))$
case this dual is the quantum group $U_q(gl(n))$ itself).

All the finite dimensional modules $V_\lambda$ over the REA are also
the modules over $\cal T^*$. The representations $\pi_\lambda$ of
the REA constructed in \cite{S} commute (as the mappings) with the
action of $\cal T^*$. We call them {\it the equivariant
representations} precisely in this sense.

At the end of the section we would like to discuss the problem
whether the category $\lhq-\Rep$ of finite dimensional equivariant
representations of the algebra $\lhq$ is "big enough".

\begin{definition} \label{def:faithful}
Let $A$ be an algebra and $A-\Rep$ be the category of its
representations. We say that this category is faithful if for any
nonzero element $a\in A$ there exists an object $V\in A-\Rep$ such
that $\pi_V(a)\not=0$ where $\pi_V$ is the representation
corresponding to $V$.
\end{definition}

For an involutive $R$-matrix (that is $R^2 = I$) the category
$\lhq-\Rep$ of equivariant representations $\pi_\lambda^z$ of the
algebra $\lhq$ is faithful. The same is true for any Hecke symmetry
$R$ which is a flat deformation of an involutive $R$-matrix (in
particular, for the $U_q(sl(n))$ $R$-matrix). This can be
established by a direct modification of the proof of the analogous
classical statement about universal enveloping algebras given in
\cite{D}. For an arbitrary Hecke symmetry we shall suppose that the
above property of the representation category takes place as a
plausible conjecture.

\section{The basic q-Newton identities}
\label{sec:5}

In this section we deal with a non-modified REA and generalize
formulae (\ref{reduc}) to this case.

First, we introduce some convenient notations. Let us denote
$\vmd(t_1,t_2,\dots,t_n)$ the following Vandermonde determinant
$$
\vmd(t_1,t_2,\dots,t_n) = \left|
\matrix{1&1&\dots &1\cr
t_1&t_2&\dots&t_n\cr
t_1^2&t_2^2&\dots&t_n^2\cr
\dots &\dots &\dots &\dots \cr
t_1^{n-1}&t_2^{n-1}&\dots&t_n^{n-1}}
\right| = \prod_{i>j}(t_i-t_j),
$$
where $t_i$ are some variables. In general, this variables can be
elements of a commutative algebra (a ring).

Besides this, we consider the elementary symmetric functions in the
variables $t_i$:
\be
e_0\equiv 1,\quad
e_k = \sum_{1\le i_1<i_2<\dots <i_k\le n}t_{i_1}t_{i_2}
\dots t_{i_k},\quad 1\le k\le n.
\label{elem-symm}
\ee
With each $e_k$ we associate the series of the following quantities
\be
e_k(\hat t_i)\stackrel{\mbox{\tiny def}}{=}{e_k}_{
\rule{0.25pt}{3.5mm}_{\;t_i=0}},\quad
e_k(\hat t_i,\hat t_j)\stackrel{\mbox{\tiny def}}{=}{e_k}_{
\rule{0.25pt}{3.5mm}_{\;t_i=0,t_j=0}},\quad 1\le i,j\le n,
\quad{\rm and\ so\ on}.
\label{hat}
\ee
It is evident, that the quantity $e_k(\hat t_i)$ is an elementary
symmetric function in the set of $(n-1)$ variables $t_j$, $j\not=i$,
etc.

Note some useful properties of the above quantities (their proof is
a simple exercise)
\begin{eqnarray}
&&e_k = e_k(\hat t_i) + t_ie_{k-1}(\hat t_i),\label{pro:1}\\
\rule{0pt}{5mm}
&&e_k(\hat t_i) - e_k(\hat t_j)  = (t_j-t_i)
e_{k-1}(\hat t_i, \hat t_j) \label{pro:2}\\
\rule{0pt}{5mm}
&& ke_k = \sum_{i=1}^{n}t_ie_{k-1}(\hat t_i).
\label{pro:3}
\end{eqnarray}
Here we assume $1\le k\le n$ and, besides, $1\le i,j\le n$ in the
first two lines.

The following lemma is easy to verify.

\begin{lemma}
\label{lem:1}
$$
\left|
\matrix{1&1&\dots &1\cr
e_1(\hat t_1)&e_1(\hat t_2)&\dots&e_1(\hat t_n)\cr
e_2(\hat t_1)&e_2(\hat t_2)&\dots&e_2(\hat t_n)\cr
\dots &\dots &\dots &\dots \cr
e_{n-1}(\hat t_1)&e_{n-1}(\hat t_2)&\dots&
e_{n-1}(\hat t_n)}
\right|=\prod_{i<j}(t_i-t_j)
\equiv \vmd(t_n,t_{n-1},\dots,t_1).
$$
\end{lemma}
\medskip

\noindent
{\bf Proof\ \ } The lemma is proved by induction in the size of the
determinant. Being based on (\ref{pro:2}), the induction proceeds in
the same way as when calculating the Vandermonde determinant.
\hfill\rule{6.5pt}{6.5pt}

\medskip
Consider now two sets of independent central elements of REA $\lq$
$$
\sigma_k(\hat L) = q^k \Tr_{R(12\dots k)}A^{(k)}\hat L_{\bar 1}\dots
\hat L_{\bar k}, \quad{\rm and}\quad s_k(\hat L) =q\,\Tr_R(\hat L^k),
\quad 1\le k\le p,
$$
were $\hat L_{\bar k}$ is defined as follows:
$$
\hat L_{\bar 1} = \hat L_1, \quad \hat L_{\bar k} =
R_{k-1}\hat L_{\overline{k-1}}\, R_{k-1}^{-1},\quad k\ge 2.
$$
We also set by definition
$$
\sigma_0(\hat L) = s_0(\hat L) = {\rm id}_{\cal L}.
$$

These two sets of central elements are connected by the Newton
identities \cite{GPS,PS}
\be
\begin{array}{rcl}
s_1 &=& \sigma_1  \\
\rule{0pt}{4mm}
-s_2+s_1\sigma_1& =&
2_qq^{-1}\,\sigma_2 \\
\rule{0pt}{4mm}
s_3-s_2\sigma_1+s_1\sigma_2 &=&
3_qq^{-2}\,\sigma_3 \\
\rule{0pt}{4mm}
\dots&\dots &\dots \\
\rule{0pt}{4mm}
(-1)^{p-1}s_p+(-1)^{p-2}s_{p-1}\sigma_1 + \dots +
s_{1}\sigma_{p-1}
&=& p_qq^{1-p}\,\sigma_p \\
\end{array}
\label{NI}
\ee
Using (\ref{NI}) one can in principle express the quantities $s_k$
(and, therefore,  ${\rm Tr}_{R}(\hat L^k)$) in terms of $\sigma_i$,
$1\le i\le k$. But the corresponding expressions are very cumbersome
and provide no advantage in working with ${\rm Tr}_{R}(\hat L^k)$.
On the other hand, there exists a useful and handy parametric
resolution of the system of Newton identities.

Namely, we shall assume the central elements $\sigma_k(\hat L)$ to
be represented in the form, analogous to (\ref{elem-symm}) (see also
(\ref{elsym}))
\be
\sigma_k(\hat L) = \sum_{1\le i_1<i_2<\dots <i_k\le p}\mu_{i_1}
\mu_{i_2} \dots \mu_{i_k},\quad 1\le k\le p.
\label{q-elem-symm}
\ee
The elements $\mu_i$, $1\le i\le p$, belong to the algebraic closure
of the center of the REA. On passing to an orbit ${\cal L}_q^\chi$,
the quantities $\mu_i$ take numerical values from the ground field.
In this case we assume that all these numbers are distinct pairwise.
So, we have a 1-generic orbit. (Hereafter, all the notions are used
by analogy with those introduced in section \ref{sec:fuzzy}.)

The above mentioned parametric resolution of (\ref{NI}) is given in
the following proposition.

\begin{proposition}
Let the central elements $\sigma_k(\hat L)$ be parametrized by
(\ref{q-elem-symm}). Then
\be
q^{-1}s_k(\hat L)\equiv {\rm Tr}_R(\hat L^k) =
q^{-p} \sum_{i=1}^p\mu_i^kd_i
\label{trl}
\ee
where
\be
d_i = \prod_{j\not=i}^p{{q\mu_i-q^{-1}\mu_j}\over{\mu_i-\mu_j}}.
\label{DoM}
\ee
\end{proposition}

\noindent
{\bf Proof\ \ } Let us denote
$$
x_i = q^{1-p}d_i =
\prod_{j\not=i}^p{{\mu_i-q^{-2}\mu_j}\over{\mu_i-\mu_j}}.
$$
Then we should prove that $s_k = \sum\mu_i^kx_i$ is a solution of
(\ref{NI}), provided that $\sigma_k$ is given by
(\ref{q-elem-symm}).

First of all, note the following representation for $x_i$:
\be
x_i = q^{(p-1)(p-2)}\,
\frac{\vmd(q^{-2}\mu_1,q^{-2}\mu_2,\dots,\mu_i,\dots,q^{-2}
\mu_p)}{\vmd(\mu_1,\mu_2,\dots,\mu_p)}.
\label{x-vmd}
\ee
It is a direct consequence of the explicit form of $x_i$.

Now we substitute the ansatz $s_k=\sum\mu_i^kx_i$ into the set of
Newton identities and prove that the corresponding system of linear
equation in the variables $x_i$ has a unique solution which
coincides with (\ref{x-vmd}).

Using (\ref{pro:1}) and (\ref{q-elem-symm}), we transform (\ref{NI})
to the following system of linear equations:
\be
\sum_{i=1}^p\mu_i\sigma_{k-1}(\hat\mu_i)x_i = k_q q^{1-k}
\,\sigma_k,
\qquad 1\le k\le p,
\label{system}
\ee
where the quantities $\sigma_k(\hat \mu_i)$ have the same meaning as
$e_k(\hat t_i)$ in (\ref{hat}).

The determinant of the system is
$$
\Delta(\mu) =
\left|
\matrix{\mu_1&\mu_2&\dots &\mu_p\cr
\mu_1\sigma_1(\hat \mu_1)&\mu_2\sigma_1(\hat \mu_2)
&\dots&\mu_p\,\sigma_1(\hat \mu_p)\cr
\mu_1\sigma_2(\hat \mu_1)&\mu_2\sigma_2(\hat \mu_2)&
\dots&\mu_p\,\sigma_2(\hat \mu_p)\cr
\dots &\dots &\dots &\dots \cr
\mu_1\sigma_{p-1}(\hat \mu_1)&\mu_2\sigma_{p-1}(\hat \mu_2)&\dots&
\mu_p\,\sigma_{p-1}(\hat \mu_p)}
\right|.
$$
Taking into account Lemma \ref{lem:1}, one can rewrite the above
determinant in an equivalent form:
\be
\Delta(\mu) = \left(\prod_{i=1}^p\mu_i\right)
\vmd(\mu_p,\mu_{p-1},\dots,\mu_1).
\label{det-sys}
\ee
Since $\Delta(\mu)\not=0$, the system (\ref{system}) has a unique
solution. To find the solution we use the Cramer's formula.
Evidently, it suffices to find the value of $x_1$ say, since the
values of other variables can be obtained from the letter one by
simple permutation of $\mu_i$.

So, we shall find $x_1$. In accordance with the Cramer's rule, it is
equal to the ratio of the following determinants:
$$
x_1 = \frac{1}{\Delta(\mu)}\,
\left|
\matrix{\sigma_1&\mu_2&\dots &\mu_p\cr
2_qq^{-1}\sigma_2&\mu_2\sigma_1(\hat \mu_2)
&\dots&\mu_p\,\sigma_1(\hat \mu_p)\cr
\displaystyle
3_qq^{-2}\sigma_3&\mu_2\sigma_2(\hat \mu_2)&
\dots&\mu_p\,\sigma_2(\hat \mu_p)\cr
\dots &\dots &\dots &\dots \cr
p_qq^{1-p}\sigma_p&\mu_2\sigma_{p-1}(\hat \mu_2)
&\dots&\mu_p\,\sigma_{p-1}(\hat \mu_p)}
\right|\equiv \frac{\Delta_1(\mu)}{\Delta(\mu)}.
$$
Let us now identically transform the numerator of the above
expression --- the determinant $\Delta_1(\mu)$.

First of all, from the first column of $\Delta_1(\mu)$ we subtract
the sum of all other columns. Using (\ref{pro:1}) and (\ref{pro:3}),
one gets for the general element of the first column
\be
k_qq^{1-k}\,\sigma_k -
\sum_{i=2}^p\mu_i\sigma_{k-1}(\hat \mu_i) = z_k\,\sigma_k(\hat
\mu_1)+z_{k-1}\,\mu_1\sigma_{k-1}(\hat \mu_1) + q^{2(1-k)}\mu_1\,
\sigma_{k-1}(\hat \mu_1).
\label{1-sum}
\ee
where for the sake of compactness we have introduced a notation
$$
z_n = n_qq^{1-n} - n.
$$
So, we find that each element of the first column of $\Delta_1(\mu)$
is the sum of several terms and therefore one can expand
$\Delta_1(\mu)$ into the sum of determinants
$$
\Delta_1(\mu) = \Delta_1'(\mu) + \Delta_1''(\mu),
$$
where the $k$-th element of the first column of $\Delta_1'(\mu)$ is
equal to $q^{2(1-k)}\mu_1 \sigma_{k-1}(\hat \mu_1)$ while the $k$-th
element of the first column of $\Delta_1''(\mu)$ contains the sum
$\eta_k(\mu)$ of two rest terms in the right hand side of
(\ref{1-sum})
$$
\eta_k(\mu)\equiv z_k\,\sigma_k(\hat \mu_1)+z_{k-1}\,\mu_1
\sigma_{k-1}(\hat \mu_1).
$$

First, consider the determinant
$$
\Delta_1''(\mu)  = \prod_{i=2}^p\mu_i
\left|
\matrix{0&1&1& \dots &1\cr
\eta_2 &\sigma_1(\hat \mu_2)&\sigma_1(\hat \mu_3)
&\dots&\sigma_1(\hat \mu_p)\cr
\eta_3 &\sigma_2(\hat \mu_2)&\sigma_2(\hat \mu_3)&
\dots&\sigma_2(\hat \mu_p)\cr
\dots &\dots &\dots &\dots \cr
\eta_p&\sigma_{p-1}(\hat \mu_2)&\sigma_{p-1}(\hat \mu_3)
&\dots&\sigma_{p-1}(\hat \mu_p)}
\right|.
$$
We shall prove that $\Delta_1''(\mu) = 0$.

Subtracting the second column consecutively from the third one,
the fourth one and so on, and taking into account (\ref{pro:2}), one
gets
$$
\Delta_1''(\mu)  = \prod_{i=2}^p\mu_i\prod_{j=3}^p(\mu_2-\mu_j)
\left|
\matrix{0&1&0& \dots &0\cr
\eta_2 &\sigma_1(\hat \mu_2)&1&\dots&1\cr
\eta_3 &\sigma_2(\hat \mu_2)&\sigma_1(\hat \mu_2,\hat \mu_3)&
\dots&\sigma_1(\hat \mu_2,\hat \mu_p)\cr
\dots &\dots &\dots &\dots &\dots\cr
\eta_p&\sigma_{p-1}(\hat \mu_2)&\sigma_{p-2}(\hat \mu_2,\hat
\mu_3)
&\dots&\sigma_{p-2}(\hat \mu_2,\hat \mu_p)}
\right|.
$$
Then we repeat this procedure subtracting the third column from each
$j$-th column with $j>3$ and so on. As a result we come to
$$
\Delta_1''(\mu)  = N(\mu)
\left|
\matrix{0 & 1 & 0 & 0 & \dots & 0 &0\cr
\eta_2 &\sigma_1(\hat \mu_2)&1& 0 &\dots&0 &0\cr
\eta_3 &\sigma_2(\hat \mu_2)&\sigma_1(\hat \mu_2,\hat \mu_3)& 1&
\dots& 0 & 0\cr
\dots &\dots &\dots &\dots &\dots &\dots &\dots\cr
\eta_{p-1}&\sigma_{p-2}(\hat \mu_2)&\sigma_{p-3}(\hat \mu_2,
\hat \mu_3)&
\sigma_{p-4}(\hat \mu_2,\hat \mu_3,\hat \mu_4)&
\dots&\sigma_{1}(\hat \mu_2,\dots)& 1\cr
\eta_p&\sigma_{p-1}(\hat \mu_2)&\sigma_{p-2}(\hat \mu_2,
\hat \mu_3) &\sigma_{p-3}(\hat \mu_2,\hat \mu_3,\hat \mu_4)
&\dots&\sigma_{2}(\hat \mu_2,\dots)& \mu_1
}
\right|
$$
with
$$
N(\mu) = \prod_{i=2}^p\mu_i\prod_{2\le j<k \le p}(\mu_j-\mu_k).
$$
The result obtained admits the further simplification. We multiply
the last column by $\mu_p$ and subtract it from the $(p-1)$-th
column, then multiply the last column by $\mu_{p-1}\mu_p$ and
subtract it from the $(p-2)$-th column, etc. Then we repeat the
procedure starting from the $(p-1)$-th column of the resulting
determinant and so on. It is not difficult to see that we end up
with the following form of $\Delta_1''(\mu)$
$$
\Delta_1''(\mu) = N(\mu)
\left|
\begin{array}{lcccccc}
0 & 1 & 0 & 0 & \dots & 0 &0\\
z_2\sigma_2(\hat\mu_1) &\mu_1&1& 0 &\dots&0 &0\\
z_3\sigma_3(\hat\mu_1)+z_2\sigma_2(\hat\mu_1) \mu_1
&0 &\mu_1& 1& \dots& 0 & 0\\
z_4\sigma_4(\hat\mu_1)+z_3\sigma_3(\hat\mu_1)\mu_1 &0 &0& \mu_1&
\dots& 0 & 0\\
\dots &\dots &\dots &\dots &\dots &\dots &\dots\\
z_{p-1}\sigma_{p-1}(\hat\mu_1)+z_{p-2}\sigma_{p-2}(\hat\mu_1)\mu_1
&0&0&0&
\dots&\mu_1& 1\\
z_{p-1}\mu_1\sigma_{p-1}(\hat\mu_1) &0&0&0&\dots&0&\mu_1
\end{array}
\right|
$$
where we have restored the explicit form of $\eta_k$ and have taken
into account that $\sigma_p(\hat\mu_1)\equiv 0$.

At last, we multiply the third column by $z_2\sigma_2(\hat\mu_1)$,
the fourth one by $z_3\sigma_3(\hat\mu_1)$, etc., and then subtract
all these columns from the first one. We get all elements of the
first column of $\Delta_1''(\mu)$ to be zero, therefore
$\Delta_1''(\mu) = 0$.

Turn now to the determinant
$$
\Delta_1'(\mu) = \prod_{i=1}^p\mu_i
\left|
\matrix{1 & 1 & 1 & \dots & 1 \cr
q^{-2}\sigma_1(\hat \mu_1) &\sigma_1(\hat \mu_2)&\sigma_1
(\hat \mu_3) &\dots&\sigma_1(\hat \mu_p)\cr
q^{-4}\sigma_2(\hat \mu_1) &\sigma_2(\hat \mu_2)&\sigma_2
(\hat \mu_3)& \dots&\sigma_2(\hat \mu_p)\cr
\dots &\dots &\dots &\dots &\dots \cr
q^{2(1-p)}\sigma_{p-1}(\hat \mu_1)&\sigma_{p-1}(\hat \mu_2)&
\sigma_{p-1}(\hat \mu_3) &\dots&\sigma_{p-1}(\hat \mu_p)}
\right|.
$$
With the same operations which were applied to $\Delta_1''(\mu)$ we
convert the determinant into the form
$$
\Delta_1'(\mu) = \prod_{i=1}^p\mu_i\prod_{2\le j<k\le p}
(\mu_j - \mu_k) \left|
\matrix{1 & 1 & 0 & \dots & 0 \cr
q^{-2}\sigma_1(\hat \mu_1) &\mu_1&1&\dots&0\cr
q^{-4}\sigma_2(\hat \mu_1) &0&\mu_1& \dots&0\cr
\dots &\dots &\dots &\dots &\dots\cr
q^{2(1-p)}\sigma_{p-1}(\hat \mu_1)&0&
0 &\dots&\mu_1}
\right|.
$$
Let us introduce new parameters
$$
\nu_1 = \mu_1,\quad \nu_i = q^{-2}\mu_i,\quad 2\le i\le p.
$$
Since the function $\sigma_k(\hat \mu_1)$ is a homogeneous
polynomial of the $k$-th order in the variables $\mu_i$, $i\ge 2$,
then
$$
q^{-2k}\sigma_k(\hat \mu_1) = \sigma_k(\hat \nu_1)
$$
and therefore
$$
\Delta_1' = \left(q^{2(p-1)}\prod_{i=1}^p\nu_i\right)q^{(p-1)(p-2)}
\left|
\matrix{1 & 1 & 1 & \dots & 1 \cr
\sigma_1(\hat \nu_1) &\sigma_1(\hat \nu_2)&\sigma_1(\hat \nu_3)
&\dots&
\sigma_1(\hat \nu_p)\cr
\sigma_2(\hat \nu_1) &\sigma_2(\hat \nu_2)&\sigma_2(\hat \nu_3)&
\dots&
\sigma_2(\hat \nu_p)\cr
\dots &\dots &\dots &\dots &\dots\cr
\sigma_{p-1}(\hat \nu_1)&\sigma_{p-1}(\hat \nu_2)&
\sigma_{p-1}(\hat \nu_3) &\dots&\sigma_{p-1}(\hat \nu_p)}
\right|.
$$
Applying Lemma \ref{lem:1} and changing the set of parameters
$\{\nu_i\}$ back to the set of $\{\mu_i\}$ we get
$$
\Delta_1(\mu) = \Delta_1'(\mu) = q^{(p-1)(p-2)}
\left(\prod_{i=1}^p\mu_i\right)
\vmd(q^{-2}\mu_p,q^{-2}\mu_{p-1},\dots,q^{-2}\mu_{2},\mu_1).
$$
Using the value (\ref{det-sys}) of the determinant $\Delta(\mu)$,
one comes to the final result
$$
x_1 = \frac{\Delta_1(\mu)}{\Delta(\mu)} =
q^{(p-1)(p-2)}\frac{\vmd(q^{-2}\mu_p,q^{-2}\mu_{p-1},\dots,\mu_1)}
{\vmd(\mu_p,\mu_{p-1},\dots,\mu_1)}
$$
which is obviously equivalent to (\ref{x-vmd}).
\hfill\rule{6.5pt}{6.5pt}

\begin{remark}
\label{DM} {\rm
In \cite{DM1} a way of quantization of semisimple (but not necessary
generic) orbit in $gl(n)^*$ was suggested and in this connection
another form of the basic $q$-Newton identity was given. Let us
briefly describe the quantization procedure from \cite{DM1} and
compare two forms of basic $q$-Newton identities. (The normalization
of the quantum trace in \cite{DM1} differs from that accepted in the
present paper.) Here we restrict ourselves  to the $\Uq$ case.

Consider the $GL(n)$ orbit ${\cal O}_M$ of an arbitrary semisimple
matrix $M\in gl(n)^*$ characterized by $r\le n$ pairwise distinct
eigenvalues $\mu_i$ with the multiplicities $m_i\geq 1$
\be
m_1+m_2+\dots +m_r=n
\label{spec-dat}
\ee
(see Remark 1).

Let $P(x)=\prod_{i=1}^r(x-\mu_i)$ be the degree $r$ minimal
polynomial of this orbit (i.e. each $\mu_i$ is a simple root of
$P(x)$). Then the quotient of the mREA over the ideal generated by
the entries of the matrix $P(L)=\prod_{i=1}^r(L-\mu_i\Id)$ and by
the elements
\be
\Tr_R L^k-\bar{\beta_k},\quad k=1,\dots ,r-1
\label{nong-orb}
\ee
with appropriate values of $\bar{\beta_k}$ is a flat deformation
(quantization) of the commutative algebra $\K({\cal O}_M)$. This
fact was established in \cite{DM1} where the exact values of
$\bar{\beta_k}$ were expressed in terms of the roots of the minimal
polynomial. However, the same values of these quantities can be
obtained from our parametric resolution of the basic $q$-Newton
identities given in (\ref{trl})--(\ref{DoM}).

For this purpose, we associate to any root $\mu=\mu_i$ of the
minimal polynomial the following {\it string} of $m_i$ quantities
\be
\nu_1=\mu,\;\;
\nu_2=q^{-2}\nu_1+q^{-1}\h,\;\;\nu_3=q^{-2}\nu_2+q^{-1}
\h,\;\;\dots,\;\;
\nu_{m_i}=q^{-2} \nu_{m_i-1}+q^{-1}\h.
\label{string}
\ee
Consider the set of all $\nu_{i_k}$ belonging to strings
(\ref{string}) and define the character $\chi:Z(\lhq)\to\K$ on
central elements $\sigma_i$ as in (\ref{q-gen}) but with eigenvalues
running over all $\nu$ from the mentioned set. Let us pass to the
corresponding NC orbit $\lhqc$ (\ref{q-orb}). Emphasize that the NC
orbit $\lhqc$ thus obtained turns out to be bigger than the result
of quantization of the Poisson pencil (see Introduction) on the
initial $GL(n)$ orbit ${\cal O}_M$. In order to get the genuine
quantum orbit we should quotient $\lhqc$ over the ideal generated by
the entries of the matrix $P(L)$ constructed from the minimal
polynomial of the orbit ${\cal O}_M$. So, assuming $|q-1|$ and $\h$
to be small enough in order to avoid any casual coincidence of the
elements from the above union of the strings, we get $n$
{\it pairwise distinct} eigenvalues $\nu_{i_k}$ of the matrix $L$
with entries considered as elements of $\lhqc$. This is well
coordinated with the empirical principle that the quantization
decreases the degeneracy.

Finally, we have got a 1-generic orbit $\lhqc$ and the quantities
$\Tr_R L^k$ can be computed via (\ref{trl})-(\ref{DoM}). However,
the multiplicities $d_i$ corresponding to extra eigenvalues (i.e.
those which are not roots of the minimal polynomial) vanish.

So, given a 1-generic NC orbit $\lhqc$, we can be sure that it is a
quantization of a classical generic orbit iff the set of eigenvalues
of the matrix $L$ corresponding to this NC orbit contains no string.
If it is not the case, we construct the minimal polynomial $P(x)$
taking the first element of each string as its simple root and
consider a two sided ideal in $\lhqc$ generated by the entries of
the matrix $P(L)$. Then, quotienting the given 1-generic NC orbit
over this ideal, we get a quantization of a semisimple but not
generic orbit whose eigenvalues are the first elements of strings
and the multiplicity of each eigenvalue is equal to the length of
the corresponding string.}
\end{remark}

\section{The higher Cayley-Hamilton and Newton identities}
\label{sec:6}

Consider a 1-generic quantum orbit $\lhqc$ defined by relations
(\ref{q-gen})--(\ref{ideal}) where $p$ is the symmetry rank of the
Hecke $R$-matrix.

In \cite{GLS2} the following conjecture was formulated.

\begin{conjecture}
On a 1-generic NC orbit $\lhqc$ the matrix $L_{(m)}$ (\ref{L-symm})
satisfies the Cayley-Hamil\-ton identity
$$
{\cal CH}_{(m)}^{\chi}(L_{(m)}) = 0,
$$
where the degree of the polynomial ${\cal CH}_{(m)}^{\chi}$ is
\be
{\rm deg\,}{\cal CH}_{(m)}^{\chi} = {m+p-1\choose m}
\label{CH-order}
\ee
and its roots $\mu_\k(m)$ are
\be
q^{m-1}\mu_\k(m) = \sum_{i=1}^p \frac{(k_i)_q}{q^{m-k_i}} \,
\mu_i + \hbar\,\xi_p(k_1,\dots,k_p),
\label{q-high}
\ee
where $\k$ is a partition (\ref{part}) of the integer $m$
$$
\k=(k_1,\dots ,k_p),\quad k_i\geq 0,\quad |\k|=k_1+\dots +k_p=m
$$
and $\xi_p(k_1,\dots,k_p)$ is the symmetric function in $k_i$ of the
form
$$
\xi_p(k_1,\dots,k_p) = \sum_{s=2}^{p} q^{k_1+k_2+\dots+k_s-m}
(k_s)_q (k_1+k_2+\dots+k_{s-1})_q.
$$
\end{conjecture}

\begin{remark}{\rm
The fact that $\xi_p$ is a symmetric function in $k_i$ can be easily
verified upon expanding all $q$-numbers in accordance with their
definition (\ref{q-num}).

Note that (\ref{q-high}) is a generalization of the analogous
formula (\ref{high}). In the classical limit $q\to 1$ the above
formula transforms into (\ref{high}) for the 1-generic fuzzy orbit
in $U(gl(p))$.}
\end{remark}

The Conjecture is justified by explicit calculations for small
values of $m$ but we still have no general proof of it. Here we
present a proof for the particular case $p=2$.

Since at $q\not=\pm 1$ the mREA $\lqh$ is isomorphic to the
nonmodified REA $\lq$, we first prove the Conjecture for $\lq$ and
then pass to the corresponding result for $\lqh$ by means of the
shift of generators.

\begin{proposition}\label{main-prop}
Let the symmetry rank $p$ of the $R$-matrix be equal to $2$. Then
the roots $\hat \omega_s$ $0\le s\le m$ of the CH polynomial for the
symmetric matrix $\hat L_{(m)}$ are given by
\be
q^{1-m} \hat \omega_s =q^{s-m}s_q\hat \mu_1+q^{-s}(m-s)_q
\hat \mu_2,\quad s=0,1,\dots,m
\label{omegas}
\ee
where $\hat \mu_i$ are the roots of the basic CH polynomial for the
matrix $\hat L_{(1)} = \|\hat l_i^{\,j}\|$ composed of the
generators of the REA $\lq$
$$
(\hat L_{(1)} - \hat \mu_{1}I_e)(\hat L_{(1)} -
\hat \mu_{2}I_e) = 0, \qquad I_e = {\rm id}_{\cal L}\otimes I.
$$
\end{proposition}
\noindent{\bf Proof\ \ }
Let us shortly outline the strategy of the proof. We shall use the
same approach as in proposition \ref{prop:3}. Namely, we prove the
claim in each finite dimensional representation $V_{(k)}$ of the
REA. This means that, for a given fixed $m$, we associate to
$\hat L_{(m)}$ a series of numerical matrices $\hat L_{(k,m)}$ of
the form
\be
\hat L_{(k,m)}
= (\bar \pi_{(k)}\otimes \pi_{(m)})(\Cas),\quad k=m,m+1,\dots
\label{def:Lkm}
\ee
and prove the claim for each value of $k\ge m$. We should consider
the values $k\ge m$ since at $k<m$ the matrix $\hat L_{(k,m)}$ can
satisfy the  CH identity of an order lower than (\ref{CH-order}).
This is a peculiarity of low dimensional representations connected
with the combinatorics of the Young diagrams.

The above matrix $\hat L_{(k,m)}$ coincides with that of a linear
operator ${\bf Cas}_{(k,m)} \in {\rm End}_r (V_{(k)})\otimes
{\rm End}_l(V_{(m)})$ which is obtained from ${\bf Cas}$
(\ref{quant-Cas}) in the full analogy with the fuzzy case
construction (\ref{cas-lr}). Recall, that at $p=2$ the spaces
$V_{(k)}$ and $V^*_{(k)}$ are isomorphic (see section
\ref{sec:rea-br}). We consider the decomposition of the
tensor product $V_{(k)}\otimes V_{(m)}$ into the direct sum of
subspaces $V_{(\nu_1,\nu_2)}$ and prove that the restriction of
${\bf Cas}_{(k,m)}$ on each subspace $V_{(\nu_1,\nu_2)}$ is a
multiple of the unit operator. The corresponding factors will be the
roots of the CH polynomial for $\hat L_{(k,m)}$ and we should verify
that they are given by (\ref{omegas}).

{\bf Step 1.} Let us fix an arbitrary integer $k\ge m$ and find the
values $\hat \mu_1(k)$ and $\hat \mu_2(k)$ of the roots of the basic
CH polynomial for
$\hat L_{(k,1)} = \|\bar\pi_{(k)}(\hat l_i^{\,j})\|$. For this
purpose we substitute the matrix $\hat L_{(k,1)}$ into identity
(\ref{CH})
$$
\hat L^2 - \sigma_1(\hat L)\hat L +\sigma_2(\hat L)\,I = 0
$$
and calculate the spectrum of the central elements $\sigma_i$. Since
$$
\sigma_1(\hat L) = q\Tr_R\hat L\quad \sigma_2(\hat L) =
\frac{q^2}{2_q}\,\Bigl(q(\Tr_R\hat L)^2 - \Tr_R\hat L^2\Bigr),
$$
then basing on the explicit form (\ref{rep:pim-r}) of the right
representation $\bar\pi_{(k)}$ one can show
$$
\sigma_1(\hat L_{(k,1)}) = (1+q^{-2k-2})\,S^{(k)}(R)\equiv
(1+q^{-2k-2})\,{\rm id}_{V_{(k)}}
$$
and
$$
\sigma_2(\hat L_{(k,1)}) = q^{-2k-2}\,S^{(k)}(R)\equiv q^{-2k-2}\,
{\rm id}_{V_{(k)}}.
$$
So, on the subspace $V_{(k)}\subset V^{\otimes k}$ the
Cayley-Hamilton identity for $\hat L_{(k,1)}$ takes the form
$$
(\hat L_{(k,1)} - I)(\hat L_{(k,1)}-q^{-2k-2}I) = 0
$$
that is
\be
\hat \mu_1(k) = 1, \quad \hat \mu_2(k) = q^{-2k-2}.
\label{mui}
\ee

{\bf Step 2.} Now we pass to the matrix $\hat L_{(k,m)}$
(\ref{def:Lkm}) treated as that of a linear operator from
${\rm End}_r (V_{(k)})\otimes {\rm End}_l(V_{(m)})$. In the case
$p=2$ the general decomposition (\ref{razl-prostr}) reduces to
$$
V_{(k)}\otimes V_{(m)} = \bigoplus_{s=0}^{m}V_{(k+s,m-s)}.
$$
The subspace $V_{(k+s,m-s)}$ can be represented as an image of the
operator ${\cal P}_s\in {\rm End}(V^{\otimes (k+m)})$
$$
{\cal P}_s = S^{(k)}_{1\dots k}S^{(m)}_{k+1 \dots k+m}
Y_{(k+s,m-s)}S^{(k)}_{1\dots k}S^{(m)}_{k+1 \dots k+m}.
$$
Here the lower indices of the $q$-symmetrizers $S^{(k)}$ and
$S^{(m)}$ indicate explicitly the numbers of factors in the tensor
product $V^{\otimes (k+m)}$ in which these symmetrizers act. The
$q$-projector $Y_{(k+s,m-s)}$ corresponds to the following Young
tableau
\be
Y_{(k+s,m-s)}\;\leftrightarrow\;
\begin{array}{|c@{\hspace{-5pt}}c|}\hline
1\;\;\;\dots \;\;\;k \;\;\; \dots\;\;\; k&+\;s\\ \hline
\multicolumn{1}{|c|}{k+s+1\;\dots\; k+m }\\ \cline{1-1}
\end{array}\equiv [1,\,\dots\,, k+s\,|\,k+s+1,\,\dots\,, k+m].
\label{Y-diag}
\ee
Here in the last equality we have introduced a more convenient
notation for an explicit enumeration of a tableau corresponding to
the two row partition $\lambda =(k+s,m-s)$, $0\le s\le m$.

Due to the fact that
$$
S^{(k)}_{1\dots k} Y_{(k+s,m-s)} = Y_{(k+s,m-s)}
S^{(k)}_{1\dots k} = Y_{(k+s,m-s)}
$$
the expression for ${\cal P}_s$ can be simplified to
\be
{\cal P}_s = S^{(m)}_{k+1 \dots k+m} Y_{(k+s,m-s)}
S^{(m)}_{k+1\dots k+m}.
\label{Ps}
\ee
The operator ${\cal P}_s$ is a projector up to a normalizing factor
$$
{\cal P}_s^2 = \gamma_s {\cal P}_s,
$$
the exact value of $\gamma_s$ is not important for us.

To prove the proposition it is suffice to show that the matrix
$\hat L_{(k,m)}$ commutes with ${\cal P}_s$ for any $0\le s\le m$
and that the following relation holds
\be
{\hat L}_{(k,m)}{\cal P}_s = \hat \omega_s{\cal P}_s =
{\cal P}_s{\hat L}_{(k,m)}.
\label{eigenv}
\ee

{\bf Step 3.} Consider the product ${\hat L}_{(k,m)}{\cal P}_s$ in
detail. On the base of definition (\ref{def:Lkm}) with the explicit
form of representations $\pi_{(m)}$ and $\bar \pi_{(k)}$ given in
propositions \ref{prop:pim} and \ref{prop:rea-r} we find the
following expression for the matrix $\hat L_{(k,m)}$
\be
q^{m-1}\hat L_{(k,m)} = \frac{m_q}{q^{2k+2}}\,S^{(k)}_{1\dots k}
S^{(m)}_{k+1\dots k+m}+\frac{\zeta m_q (k+1)_q}{q^{k+1}}\,
S^{(m)}_{k+1\dots k+m} S^{(k+1)}_{1\dots k+1}
S^{(m)}_{k+1\dots k+m}.
\label{Lkm-expl}
\ee
where we have extracted the overall factor $q^{m-1}$ to simplify the
subsequent calculations. In deriving this formula one should use the
relation
$$
S^{(k)}_{1\dots k}A^{(2)}_{kk+1}S^{(k)}_{1\dots k} =
\frac{(k+1)_q}{2_qk_q} \,(S^{(k)}_{1\dots k}
I_{\raisebox{-1.75pt}{\scriptsize $k+1$}}-S^{(k+1)}_{1\dots k+1}).
$$
Now we are to multiply the above expression for $\hat L_{(k,m)}$ on
the matrix ${\cal P}_s$ from the right. The multiplication of the
first summand  in (\ref{Lkm-expl}) results in
$q^{-2k-2}m_q{\cal P}_s$ and the main difficulty concentrates in the
second term
\be
\frac{\zeta m_q (k+1)_q}{q^{k+1}}\, S^{(m)}_{k+1\dots k+m}
S^{(k+1)}_{1\dots k+1}S^{(m)}_{k+1\dots k+m} Y_{(k+s,m-s)}
S^{(m)}_{k+1\dots k+m} \equiv \frac{\zeta m_q (k+1)_q}{q^{k+1}}
\Omega_s
\label{Om}
\ee
were the left hand side is a definition of the symbol $\Omega_s$. At
$s=m$ this is obviously the multiple of ${\cal P}_m\equiv S^{(k+m)}$
therefore below we shall suppose $0\le s\le m-1$.

To prove that $\Omega_s$ is  a multiple of ${\cal P}_s$ one should
somehow get rid of the $q$-symmetrizer $S^{(k+1)}_{1\dots k+1}$ in
its expression. For this purpose we decompose the projector
$S^{(m)}$ standing between $S^{(k+1)}$ and $Y_{(k+s,m-s)}$ into the
product of the factors (\ref{q-symm-S})
$$
S^{(m)}_{k+1\dots k+m} = \frac{(m-1)_q}{m_q}\,
S^{(m-1)}_{k+2\dots k+m} \Bigl(\frac{q^{1-m}}{(m-1)_q}\,I +
R_{k+1}\Bigr)S^{(m-1)}_{k+2\dots k+m}.
$$
The first $S^{(m-1)}$ in the right hand side of the above relation
commutes with $S^{(k+1)}$ and can be absorbed into the most left
$S^{(m)}$ in (\ref{Om})
$$
S^{(m)}_{k+1\dots k+m}S^{(m-1)}_{k+2\dots k+m} =
S^{(m)}_{k+1\dots k+m}.
$$
As for the second $S^{(m-1)}$, we shall continue its decomposition
in the same way until descending to $S^{(m-s)}$ which is absorbed by
the projector $Y_{(k+s,m-s)}$. So, we come to the result
$$
\Omega_s = \frac{(m-s)_q}{m_q} \,S^{(m)}S^{(k+1)}
\Bigl(\frac{q^{1-m}}{(m-1)_q}\,I +R_{k+1}\Bigr)\dots
\Bigl(\frac{q^{s-m}}{(m-s)_q}\,I+R_{k+s}\Bigr)Y_{(k+s,m-s)}
S^{(m)}.
$$

The next step is to draw the factors of the type $(zI+R)$ through
$Y_{(k+s,m-s)}$ and then absorb them into the right
$S^{(m)}_{k+1\dots k+m}$ on the base of the following property of
the $q$-symmetrizer
$$
R_{k+i}\,S^{(m)}_{k+1\dots k+m} = q\, S^{(m)}_{k+1\dots k+m}\qquad
1\le i\le m-1.
$$
The commutation of the linear in $R$-matrix terms with
$Y_{(k+s,m-s)}$ is done on the base of the identity \cite{OgP}
\be
\Bigl(\frac{q^{-i}}{i_q}\,I +R_i\Bigr)\;
\begin{array}{|cc|}\hline
1\;\;\;\dots\;\;\; i&\dots\\ \hline
\multicolumn{1}{|c|}{i+1\;\;\dots}
\\ \cline{1-1}
\end{array}
=
\begin{array}{|cc|}\hline
1\;\;\;\dots& i+1\;\;\dots\\ \hline
\multicolumn{1}{|l|}{i\;\;\dots}
\\ \cline{1-1}
\end{array}
\;\Bigl(R_i - \frac{q^{i}}{i_q}\,I \Bigr) \label{link}
\ee
On applying such like formulae $s$ times we get the following result
for $\Omega_s$
\be
\Omega_s = \frac{(m-s)_q}{(k+s)_q}\sum_{r=1}^s
\frac{(k-m+2r)_q}{(m-r)_q(m-r+1)_q}\,S^{(m)}Y_{(k+s,m-s)}^{[r]}
S^{(m)}, \quad 0\le s\le m-1.
\label{Om-prom}
\ee
Here the $q$-projector $Y_{(k+s,m-s)}^{[r]}$ corresponds to the
following Young tableau
$$
Y_{(k+s,m-s)}^{[r]}\leftrightarrow [1,\,\dots\,,k+r,k+r+2,\,
\dots\,,k+s+1\,|\,k+r+1,k+s+2,\,\dots\,,k+m].
$$
The unwanted $q$-symmetrizer $S^{(k+1)}$ has been absorbed into
$Y^{[r]}$
$$
S^{(k+1)}_{1\dots k+1}Y^{[r]}_{(k+s,m-s)} = Y^{[r]}_{(k+s,m-s)},
\quad 1\le r\le s.
$$

And at last, we transform all $Y^{[r]}_{(k+s,m-s)}$ in
(\ref{Om-prom}) back to $Y_{(k+s,m-s)}$. This can be done by the
multiple successive application of the following consequence of
(\ref{link})
$$
\begin{array}{|cc|}\hline
1\;\;\;\dots& i+1\;\;\dots\\ \hline
\multicolumn{1}{|l|}{i\;\;\dots}
\\ \cline{1-1}
\end{array} =
\frac{i_q^2}{(i-1)_q(i+1)_q}\, \Bigl(\frac{q^{-i}}{i_q}\,I
+R_i\Bigr)\;
\begin{array}{|cc|}\hline
1\;\;\;\dots\;\;\; i&\dots\\ \hline
\multicolumn{1}{|c|}{i+1\;\;\dots}
\\ \cline{1-1}
\end{array}
\;\Bigl(\frac{q^{-i}}{i_q}\,I + R_i\Bigr).
$$
The brackets with $R$-matrices appearing in this way are absorbed
into $S^{(m)}$ which gives rise to the accumulation of an overall
numerical factor. The final result for $\Omega_s$ reads
$$
\Omega_s = \beta_s \,S^{(m)} Y_{(k+s,m-s)}S^{(m)}, \quad
0 \le s \le m-1,
$$
where
$$
\beta_s = (m-s)_q(k+s+1)_q\sum_{r=1}^s
\frac{(k-m+2r)_q}{(k+r)_q(k+r+1)_q(m-r)_q(m-r+1)_q}.
$$
The sum in this relation can be easily calculated by induction in
$s$ and one gets
$$
\beta_s = \frac{s_q(k+s+1-m)_q}{m_q(k+1)_q}.
$$

Now, gathering together the results from the both terms in the right
hand side of (\ref{Lkm-expl}), we find
$$
q^{m-1}\hat L_{(k,m)}{\cal P}_s = \Bigl(\frac{m_q}{q^{2k+2}}
+\frac{\zeta}{q^{k+1}}\, m_q(k+1)_q\beta_s\Bigr) {\cal P}_s =
\Bigl(\frac{m_q}{q^{2k+2}} +
\frac{\zeta}{q^{k+1}}s_q(k+s+1-m)_q\Bigr) {\cal P}_s.
$$
Moreover, since the expressions for $\hat L_{(k,m)}$
(\ref{Lkm-expl}) and ${\cal P}_s$ (\ref{Ps}) are symmetric with
respect to the $q$-projectors, the same result can be obtained for
${\cal P}_s\hat L_{(k,m)}$ that is ${\cal P}_s$ and $\hat L_{(k,m)}$
commute.

The last step  is to show that the coefficient in the right hand
side of the above expression for $\hat L_{(k,m)}{\cal P}_s $ is
equal to $\hat \omega_s$. It is a matter of a short straightforward
calculation to verify the identity
$$
\frac{m_q}{q^{2k+2}} + \frac{\zeta}{q^{k+1}}s_q(k+s+1-m)_q
\equiv \frac{s_q}{q^{m-s}}+\frac{(m-s)_q}{q^{2k+2+s}}.
$$
Taking into account the values of the roots $\hat \mu_i(k)$
(\ref{mui}) we come to the desired result
$$
\frac{s_q}{q^{m-s}}+\frac{(m-s)_q}{q^{2k+2+s}} = q^{s-m}{s_q}\hat
\mu_1(k)+q^{-s}(m-s)_q \hat \mu_2(k) = q^{m-1}\hat \omega_s
$$
and therefore
$$
\hat L_{(k,m)}{\cal P}_s = \hat \omega_s {\cal P}_s.
$$
As was mentioned above, the same result is valid for
${\cal P}_s \hat L_{(k,m)}$. \hfill\rule{6.5pt}{6.5pt}
\smallskip

The roots of the CH polynomial for the symmetric matrix $L_{(m)}$
composed of the generators $l_i^{\,j}$ of the mREA can be found as a
simple corollary of the proposition \ref{main-prop}.

\begin{corollary}
Let the symmetry rank $p$ of the $R$-matrix be equal to $2$. Then
the roots $\omega_s$ $0\le s\le m$ of the CH polynomial for the
symmetric matrix $L_{(m)}$ are given by
\be
q^{1-m}  \omega_s =q^{s-m}s_q \mu_1+q^{-s}(m-s)_q\mu_2 +
\h s_q(m-s)_q,\quad s=0,1,\dots,m
\label{omgas}
\ee
where $\mu_i$ are the roots of the basic CH polynomial for the
matrix $L_{(1)} = \|l_i^{\,j}\|$ composed of the generators of the
mREA $\lqh$
$$
( L_{(1)} -  \mu_{1}I_e)(L_{(1)} -  \mu_{2}I_e) = 0,
\qquad I_e = {\rm id}_{\cal L}\otimes I.
$$
\end{corollary}
\noindent {\bf Proof\ \ }
Consider in detail how the shift of generators presented in
(\ref{shift}) affects the CH identity. If the matrix
$\hat L = \|\hat l_i^{\,j}\|$ of the REA generators satisfies the CH
identity
$$
(\hat L - \hat \mu_1 I_e)(\hat L - \hat \mu_2 I_e) = 0,
$$
then the matrix $L$ of the mREA generators satisfies the same
identity but with $\mu_i = \hat\mu_i + \hbar\zeta^{-1}$. This is a
trivial consequence of (\ref{shift}).

As for the higher order CH identity for symmetrical matrix
$L_{(m)}$ (\ref{L-symm}) the modification is as follows. Suppose, we
know the CH identity for $\hat L_{(m)}$
$$
\prod_{s=0}^{m}(\hat L_{(m)}- \hat\omega_s I_e^{(m)}) = 0,\quad
q^{m-1}\hat\omega_s =  q^{s-m}s_q\hat \mu_1+q^{-s}(m-s)_q\hat \mu_2,
$$
where $I_e^{(m)} = {\rm id}_{\cal L}\otimes I_{V_{(m)}}$. To pass to
the mREA case we should take into account the connection of the
matrices $\hat L_{(m)}$ and $L_{(m)}$
$$
\hat L_{(m)} = L_{(m)} - \frac{\hbar}{\zeta}\,q^{1-m}
m_q\,I_e^{(m)}.
$$
This relation follows from (\ref{shift}), (\ref{L-symm}) and
the explicit form of the representation $\pi_{(m)}$ given in
proposition \ref{prop:pim}. Besides, one should take into account
the shift from $\hat \mu_i$ to $\mu_i$ described above. Therefore,
we come to the following result
\begin{eqnarray*}
\hat L_{(m)} - \hat \omega_sI_e^{(m)}&=&L_{(m)} - q^{1-m}
\Bigl(\frac{\hbar}{\zeta}\,m_q + q^{s-m}s_q(\mu_1 -
\frac{\hbar}{\zeta}) +q^{-s}(m-s)_q(\mu_2 -
\frac{\hbar} {\zeta}) \Bigr)\,I_e^{(m)}\\
&=& L_{(m)} - q^{1-m}(q^{s-m}s_q\,\mu_1+q^{-s}(m-s)_q\,\mu_2 +
\hbar s_q(m-s)_q)I_e^{(m)} \\
&=& L_{(m)} - \omega_s I_e^{(m)}.
\end{eqnarray*}
This completes the proof.\hfill\rule{6.5pt}{6.5pt}
\smallskip

To sum up, we conclude that on any 1-generic NC orbit $\lhqc$
(\ref{ideal}) the extended matrix $L_{(m)}$ satisfies the
$(m+1)$-th order polynomial identity with the roots (\ref{omgas})
($p=2$). Fixing the finite dimensional representation
$\bar\pi_{(k)}$ affects only the particular form of the roots
$\mu_i$ of the basic CH polynomial for the matrix $L$.

In full analogy with the constructions of section \ref{sec:fuzzy},
we can introduce idempotents $e_\k(m)$ on any $m$-generic NC orbit.
Moreover, the following proposition holds true.

\begin{proposition}
\label{prop:21}
For any m-generic NC orbit $\lhqc$ defined by
(\ref{q-gen})--(\ref{ideal}) we have
\be
\Tr_R \lm^s=q^{-p}\sum_{|\k|=m}\mu_\k(m)^s\,d_\k(m)
\label{lm-tr}
\ee
with $\mu_\k$ introduced in (\ref{q-high}). If the symmetry rank $p$
of $R$ is equal to 2, then
\be
d_\k(m)=\prod_{1\le i< j\le
p}{\frac{q^{k_i-k_j}\mu_i-q^{k_j-k_i}\mu_j-\h\,(k_i-k_j)_q}
{\mu_i-\mu_j}}.
\label{mulp}
\ee
If the symmetry rank of $R$ $p>2$ then the above expression for
$d_\k(m)$ is valid provided that the category of finite dimensional
equivariant representations of the corresponding mREA \lqh{} is
faithful (see definition \ref{def:faithful}).
\end{proposition}

{\bf Proof\ \ } The proof of (\ref{lm-tr}) is straightforward. 
As for the proof of (\ref{mulp}), it is sufficient to establish this
formula for $\h=0$. Similarly to the proof of proposition
\ref{prop:6} we consider $\opi_\lambda(\lm)$, compute the 
corresponding multiplicities $d_\k(\la,m) = \opi_\lambda(d_\k(m))$
and prove relation (\ref{mulp}) with $\mu_i = \mu_i(\lambda)$.
Then, due to the faithfulness of the representation category, we
conclude that (\ref{mulp}) takes place at the level of the algebra
itself.

So, assuming the representation $\opi_\lambda$ to be $m$-admissible
(see definition \ref{def:m-ad}) we find
\be
d_\k(\la, m)={{\dim_q (V_{\la}^*\ot\vm)_\k}\over{\dim_q
V_{\la}^*}}
\label{raz}
\ee
where $(V_{\la}^*\ot\vm)_\k$ is the irreducible component with the
label $\k$ (\ref{part}) in the tensor product $V_{\la}^*\ot\vm$ and
the $q$-dimension is defined in (\ref{def:q-dim}). This formula is
the $q$-analog of (\ref{dk}) and can be obtained by the same method.

Now let us take into account that the decomposition rules for the
tensor product of the \lqh{} modules $V_\lambda$ are the same as
those for the $U(gl(p))$ modules, $p$ being the symmetry rank of $R$
(see the beginning of section \ref{sec:rea-br}). Bearing in mind the
property {\bf C4)} of the Hecke symmetry $R$ (section
\ref{sec:rea-def}), we conclude that $V^*_\lambda$ is isomorphic to
$V_{\hat{\lambda^{\!*}}}$ as a vector space where 
$\hat{\lambda^{\!*}}$ is defined in (\ref{lambda-star}) and
(\ref{lambda-hat}). And secondly, 
$(V_{\la}^*\ot V_{(m)})_\k\cong V_{\hat{\lambda^{\!*}}+\k}$ (as
vector spaces). Since we are interested in calculating the 
$q$-dimensions which are the same for isomorphic spaces, then in all
formulae we can change the spaces $V^*_\lambda$ and 
$(V_{\la}^*\ot V_{(m)})_\k$ for the corresponding isomorphic spaces.

On the other hand, by virtue of (\ref{DoM}) and the isomorphism
mentioned above, we have in the case $m=1$
\be
d_j(\la,1)=\prod_{i\not=j}{{q\mu_j(\la)-q^{-1}q\mu_i(\la)}
\over{\mu_j(\la) - \mu_i(\la)}}= {{\dim_q (V_{\hat{\la^{\!*}}}
\ot V)_j} \over{\dim_q V_{\hat{\la^{\!*}}}}},\quad  j=1,\dots ,p.
\label{d-j}
\ee
The explicit form of the $q$-dimension reads
\be
\dim_q \vl=s_\la(q^{p-1}, q^{p-3},\dots ,q^{-p+3}, q^{-p+1}) =
\prod_{1\le i<j\le p}{{(\la_i-\la_j-i+j)_q} \over({j-i})_q}
\label{dva}
\ee
where $s_\la$ is the Schur function corresponding to the partition
$\la$. The first equality in (\ref{dva}) was established in
\cite{GLS1}. The second one can be proved by a straightforward
calculations. Its equivalent form can be also found in \cite{Ma}.

So, on taking into account (\ref{d-j}) and (\ref{dva}) we come to
the system of $p$ equations
\be
\prod_{i\not=j}{{q\mu_j(\la)- \qq \mu_i(\la)} \over{\mu_j(\la)-
\mu_i(\la)}} = \prod_{i\not=j}{{(\la_{p-i+1}-\la_{p-j+1} +
i-j+1)_q} \over {(\la_{p-i+1}-\la_{p-j+1} +i-j)_q}},\quad
1\le j\le p.
\label{razm}
\ee
Let us find $\mu_i(\la)$ from this system. It is easy to see that
there exist a solution of the form
\be
\mu_i(\la)=\eta(\la) q^{-2(\la_{p-i+1}+i)},
\label{mu-sol}
\ee
where $\eta(\lambda)$ is an arbitrary nonzero multiplier. Indeed,
given such $\mu_i(\la)$ and taking into account definition
(\ref{q-num}) of the $q$-numbers, we get
$$
{\frac{q\mu_j(\la)-q^{-1}\mu_i(\la)}
{\mu_j(\la)-\mu_i(\la)}} =
{\frac{(\la_{p+1-i}-\la_{p+1-j}+i-j+1)_q}
{(\la_{p+1-i}-\la_{p+1-j}+i-j)_q}}.
$$
At $p=2$ the above solution (\ref{mu-sol}) is unique since the system
(\ref{razm}) is actually linear in $\mu_i(\lambda)$.

\begin{conjecture}
At an arbitrary value of $p$ the solution (\ref{mu-sol})
of the system (\ref{razm}) is unique.
\end{conjecture}

Up to this Conjecture we can extend relation (\ref{mulp}) for an
arbitrary $p\ge 2$ (at $p=2$ it is true rigorously). Indeed, from
(\ref{raz}) it follows that
$$
d_\k(\la, m)=\frac{\dim_qV_{\hat{\la^{\!*}}+\k}}
{\dim_qV_{\hat{\la^{\!*}}}}.
$$
On the other hand, taking $\mu_i$ as in (\ref{mu-sol}) we find
$$
{\frac{q^{k_i-k_j}\mu_i(\la)-q^{k_j-k_i}\mu_j(\la)}
{\mu_i(\la)-\mu_j(\la)}} =
{\frac{(\la_{p+1-j}-\la_{p+1-i}+k_i-k_j-i+j)_q}{(
\la_{p+1-j}-\la_{p+1-i}-i+j)_q}}.
$$
This relation together with above form of $d_\k(\la, m)$ and 
(\ref{dva}) entails result (\ref{mulp}) for $\h=0$. A passage
to the general case can be performed by the shift of generators
(\ref{shift}). 

To complete the proof, we point out, that since the category of
finite dimensional equivariant representation of \lqh{} is faithful,
then relation (\ref{mulp}), being valid in all representations, must
take place at the level of the algebra itself.\hfill
\rule{6.5pt}{6.5pt}

\medskip

The above relation (\ref{mu-sol}) shows which values of $\mu_i$
(up to a normalizing factor) correspond to a finite dimensional
representation $\pi_\lambda$ of the NC orbit ${\cal L}_q^\chi$. To
fix the factor $\eta(\lambda)$ one should consider the value of
$\pi_\lambda(\Tr_RL)$ for the mREA \lqh{} and use the connection
(\ref{shift}). The following corollary of proposition \ref{prop:21}
holds true.

\begin{corollary}
Let ${\cal L}_{q,1}^\chi$ be the NC orbit defined by
(\ref{q-gen})--(\ref{ideal}). In the finite dimensional
representation $\pi_\lambda$ of the mREA ${\cal L}_{q,1}$ 
the eigenvalues $\bar \mu_i(\lambda)$ of the matrix $L$ of the mREA 
generators on the orbit ${\cal L}_{q,1}^\chi$ are as follows
\be
\bar\mu_i(\lambda) = \frac{(\lambda_{p-i+1}+i-1)_q}{
q^{\lambda_{p-i+1}+i-1}}.
\label{q-mu-mrea}
\ee
Here the eigenvalues $\bar\mu_i$ are the roots of the CH polynomial 
for $L$
$$
\prod_{i=1}^p(L-\bar\mu_i\,I) \equiv 0.
$$
\end{corollary}

\noindent
{\bf Proof\ \ } 
As was found in \cite{S}, the central element ${\rm Tr}_R L$ has the
following form in the representation $\pi_\lambda$
$$
\pi_\lambda({\rm Tr}_R L) = \chi_1{\rm id}_{V_\lambda}, \quad
\chi_1 = q^{-2p}\sum_{r=1}^pq^{2r-1-\lambda_r} (\lambda_r)_q.
$$
On the other hand, by virtue of (\ref{trl}), (\ref{q-elem-symm}) and
(\ref{NI}) the $R$-trace of the matrix $\hat L$ of the generators of
non-modified REA can be presented as follows
$$
\pi_\lambda({\rm Tr}_R\hat L) = \frac{1}{q}\sum_{i=1}^p\mu_i\,
{\rm id}_{V_\lambda}
$$
were $\mu_i$ are given by (\ref{mu-sol}). Taking into account the
above expressions for the traces and the fact that $L$ and $\hat L$
are connected by means of shift (\ref{shift}), we find
$$
q^{-1}\eta(\lambda)\sum_{r=1}^pq^{-2(\lambda_{p-i+1}+i)} = 
q^{-2p}\sum_{r=1}^pq^{2r-1-\lambda_r}(\lambda_r)_q -
\frac{p_q}{q^p(q-q^{-1})}.
$$
This entails 
$$
\eta(\lambda) = -  \frac{q^2}{(q-q^{-1})}.
$$
And at last, using the connection 
$\bar\mu_i = \mu_i + (q-q^{-1})^{-1}$ and values (\ref{mu-sol}), we
come to (\ref{q-mu-mrea}) which is the $q$-generalization of the
classical result (\ref{mu}).\hfill \rule{6.5pt}{6.5pt}

\section{q-Euler characteristic and group $\KK(\lhqc)$}
\label{sec:q-eiler}

In what follows we shall consider the multiplicities $d_\k(m)$ on
the set of all  generic NC orbits as functions in $\mu_i$. As was
pointed out in remark \ref{DM}, there exist orbits on which some of
the basic multiplicities $d_i$ can vanish. But even if all $d_i$ are
nonzero, some of the higher multiplicities $d_\k(m)$ can vanish as
well. We shall restrict ourselves to the set of generic NC orbits
such that $d_\k(m)\not= 0$ for all $m=1,2,\dots $ and all for all
partitions $\k\vdash m$. Such NC orbits will be called {\em strictly
generic}. However, in considering the multiplicities as functions in
$\mu_i$ it does not matter whether we exclude a low-dimensional
subset of values of $\mu_i$ (corresponding to non-strictly generic
orbits).

Let $M_\k=M_\k(m)$ be the left (for the definiteness) projective
module corresponding to the idempotent $e_\k(m)$ defined on a NC
orbit $\lhqc$ in full analogy with (\ref{high-idemp}).

For a strictly generic orbit we consider the set of projective
modules $M_\k$ with the assignment
\be
M_\k\mapsto \Tr_R \,e_\k(m).
\label{char}
\ee

In the classical limit $q=1$, $\h=0$ of the $\Uq$ case this
assignment is out of interest (its value is identically equal to 1
since the modules $M_\k$ correspond to line bundles). But for
generic $q$ and $\h$ it is not so. Let us show that the
characteristic (\ref{char}) is close to the Euler characteristic
of a line bundle over the flag variety $\Fl(\C^n)$.

It is known that the set of $SU(n)$-equivariant line bundles over
the flag variety is in the one-to-one correspondence with the set of
holomorphic one-dimensional representations of the torus
$T\subset SU(n)$. Therefore, the bundles can be labelled by the
vectors $\k=(k_1,\dots ,k_n), k_i \in \Z$ since any such a
representation is of the form
$T\ni (t_1,\dots ,t_n) \mapsto \prod t_i^{k_i}$. Two vectors
$\k=(k_1,\dots ,k_n)$ and $\k'= (k'_1,\dots ,k'_n)$ give rise to the
same representation and hence to the same line bundle iff
$k_i=k'_i+a$ with an integer $a$ (shortly iff $\k=\k'+a$) since
$\prod t_i=1$. So, applying (if necessary) a shift by an integer we
can assume that the labels $\k$ of these line bundles obey the
condition $k_i\geq 0$. It can be easily shown (cf. \cite{L}) that
the Euler characteristic of the line bundle corresponding to
$\k=(k_1,\dots ,k_n)$ equals
\be
\prod_{i<j}{(k_i-k_j+i-j)\over (i-j)}.
\label{Ec}
\ee

Now let us pass to the NC orbits. Recall the notion of the
$q$-index introduced in \cite{GLS2}. It is defined as the following
paring
$$
\langle M_\k, \opi_\lambda\rangle=\mbox{\bf \sf Tr}\,
\opi_\lambda(e_\k(m))
$$
of a projective $\lhqc$ module $M_\k$ and a finite dimensional
equivariant representation $\opi_\lambda$. (Note that the character
$\chi$ depends on the representation $\opi_\lambda$.) Here
$\mbox{\bf \sf Tr}$ stands for the categorical trace applied to the
space $\End(\vl)\ot \End(\vm)$. (In contrast with the quantities
$\Tr_R \Lm^k$ considered above where $\Tr_R$ is applied to the
second factor only.)

As follows from constructions of the previous section, for any
$m$-admissible $\lambda$ we have
\be
\langle M_\k(m), \opi_\la \rangle = \prod_{i<j} {(\la_i-\la_j+k_i-
k_j+i-j)_q \over (i-j)_q}.
\label{Hir}
\ee

Putting $\la_i=0$ in (\ref{Hir}) we get the quantity\footnote{Since
$\la=(0,0,\dots ,0)$ is not an $m$-admissible partition we cannot
directly apply the formula (\ref{Hir}) to the corresponding
representation. However, we consider the extension of (\ref{Hir}) to
this point.}
$$
\chi_q(M_\k)=\prod_{i<j}{(k_i-k_j+i-j)_q\over (i-j)_q}
$$
which is a $q$-analog of (\ref{Ec}).  We call it {\em the $q$-Euler
characteristic}.

The characteristic $\chi_q(M_\k)$ does not depend on a concrete form
of the initial braiding $R$. Note that $\chi_q(M_\k)$ is introduced
without any holomorphic structure which is usually employed in the
definition of the Euler characteristic of a line bundle since the
operator $\bar \partial$ is essentially involved in the definition
(cf., e.g. \cite{H}).

In this connection we  recall that a generic orbit in $su(n)^*$
(treated as a real algebraic variety) can be equipped with
different complex structures. These structures are labelled by
elements of the Weyl group and the Euler characteristic of a given
line bundle essentially depends on the choice of such an element. In
our setting the quantity $\chi_q(M_\k)$ depends on the ordering of
the eigenvalues $\mu_i$. These orderings are also labelled by
elements of the Weyl group.

\begin{remark} {\rm Considering a semisimple orbit ${\cal O}$ in
$su(n)^*$ to be a real algebraic variety, we treat its coordinate
ring as an $\R$-algebra. The algebra arising from the quantization
of the Kirillov bracket on such a real variety can be realized as a
quotient of the enveloping algebra $U(su(n)_\h)$ (cf. \cite{GS2}).
It is also treated as an $\R$-algebra. However, for the
corresponding quotients of the mREA it is not possible to get rid of
complex numbers (cf. \cite{DGH} where an example of the quantum
sphere is considered). That is why we treat the above quotient as a
$\C$-algebra and consider it as the quantization of the
complexification of the orbit ${\cal O}$. }
\end{remark}

Considering the $q$-Euler characteristic as a function in  $\mu_i$
we can see that
$$
\chi_q(M_\k)=\chi_q(M_{\k'})\quad \Leftrightarrow \quad \exists
\,a\in\Z:\;\k=\k'+a.
$$
Being motivated by this property, we introduce the following
equivalence relation on the set of our projective modules
$$
M_\k\sim M_{\k'}\quad \Leftrightarrow\quad \exists
\,a\in\Z:\; \k=\k'+a.
$$
The equivalence class of the module $ M_\k$ will be denoted
$[M_\k]$.

\begin{remark}{\rm In NC geometry one usually employs the following
type of equivalence of projective modules (cf. \cite{Ro}). Two
$\A$-modules are called equivalent if the corresponding idempotents
(being extended by 0 if necessary) are similar. This equivalence is
compatible with the usual trace. Namely, for two such idempotents we
have $\Tr\,\pi(e_1)=\Tr\,\pi(e_2)$ where $\pi$ is any finite
dimensional representation of the algebra $\A$. However, such an
equivalence is not compatible with the categorical trace $\Tr_R$.
The matter is that the categorical trace is not invariant with
respect to a similarity transformation $M\to P^{-1}\, M\, P$ since
the matrices $B$ and $C$ (and their higher extensions) entering
definition (\ref{q-sled}) of the categorical trace are not invariant
under this map (if $R$ differs from  the usual flip).}
\end{remark}

Now, taking the flag variety as a pattern, we define a product on
the set of the projective modules $M_\k$:
$$
M_\k\cdot M_{\k'}=M_{\k+\k'}.
$$
It should be pointed out that the above relation does not mean any
product of {\it elements} of the modules in question.

It is evident that thus defined product of modules factorizes to the
set of equivalence classes. This set equipped with the product
becomes a group. The role of the unity is played by the class of the
trivial module. This group is an analog of the Picard group of the
flag variety.

In a similar way we can introduce an analog of G-equivariant
K-theory. Usually the algebra $K_G(pt)$ of a point is defined as an
algebra of all virtual representations of the group $G$. Let $[V]$
be the equivalence class of a finite dimensional module $V$ over the
algebra $\lhq$. Consider the $\Zq$-algebra additively  generated by
all these classes where $\Zq$ is in turn a $\Z$-algebra generated by
all $q$-numbers $m_q$, $m\in \Z$ (\ref{q-num}). As the sum of
equivalence classes we take $[U]+[V] = [U\oplus V]$ and as the
product --- $[U]\cdot [V]=[U\ot V]$. The quotient of this algebra
over the subalgebra generated by
$$
[V]-\dim_q V\,[1]
$$
will be denoted $\KK(pt)$. Here $[1]$ is the class corresponding to
the trivial object $V=\K$. Thus, the algebra $\KK(pt)$ is smaller
than $K_G(pt)$ and it is somewhat similar to the group $K(pt)$ in
the non-equivariant K-theory.

Turning to the NC orbit $\lhqc$ we treat the class $[V]$ as that of
the free $\lhqc$ module $\lhqc\ot V$. The corresponding idempotent
is the identity operator $\Id_V$ and its trace equals $\dim_q V$.

Note,  that a direct sum of two projective modules is a projective
module as well. The corresponding idempotent is a direct sum of
idempotents related to the summands.  However, the modules $M_\k(m)$
and $M_{\k'}(m)$ with the same $m$ admit the usual sum. The
corresponding idempotents are added as equal size matrices. So, we
get
$$
\sum_{|\k|=m} \,e_\k(m)=\Id_{\vm}.
$$
Regardless of the addition type we associate to the sum of two or
more modules the sum of their equivalence classes.

Now, we consider the vector space spanned by the classes $[M_\k]$
and $[1]$ with coefficients from  $\Zq$. Its quotient over subspace
generated by
$$
\sum_{|k|=m} [M_\k]-\dim_q\vm [1]
$$
will be denoted $\KK(\lhqc)$. Moreover, defining $[M_\k]\cdot [1]$
to be $[M_\k]$ we equip $\KK(\lhqc)$ with a $\Z_q$-algebra
structure. Thus, we have  $[M_\k]\cdot [V]=\dim_q V\,[M_\k]$ (i.e.
we define the class of the tensor product $M_\k\ot V$ to be
$\dim_q V\,[M_\k]$).

Proofs of the following two propositions are easy and left to the
reader.

\begin{proposition}
The algebra $\KK(\lhqc)$ is generated by the unity $[1]$ and $p$
generators $x_i=[M_i],\,\,1\le i\le p$ subject to the following
relations
\be
\sum_{1\le i_1< i_2<\dots< i_k\le p} x_{i_1} \cdot \dots \cdot
x_{i_k} = \dim_q V_{(1^k)}\,[1]= {p_q! \over {k_q!(p-k)_q!}} [1].
\label{Weyl}
\ee
Here as usual $p_q!=1_q\,2_q\dots p_q$.
\end{proposition}

Canceling $[1]$ in this formula we get the relations close to those
arising in quantum co\-ho\-mo\-lo\-gy theory in a completely
different context.

\begin{example}
Let $p=3$ and $x,\, y,\, z$ be the generators of the algebra
$\KK(\lhqc)$. Then they satisfy the system
$$
x+y+z=3_q,\quad x\cdot
y+x \cdot z + y \cdot z = 3_q,\quad x\cdot y \cdot z =1.
$$
\end{example}

\begin{proposition}
The $q$-Euler characteristic being extended at $[1]$ by
$\chi_q([1])=1$ gives rise to a linear map
$$
\chi_q:\; \KK(\lhqc)\to \Z_q.
$$
\end{proposition}

\end{document}